\documentclass[final,reqno]{siamltex}
\pdfoutput=1
\usepackage{mathtools}

\usepackage[utf8]{inputenc}
\usepackage{amsfonts,amssymb, bm}
\usepackage{graphicx,color,epstopdf}
\usepackage[colorlinks,
            linkcolor=red,
            anchorcolor=blue,
            citecolor=green
            ]{hyperref}
\usepackage[inline]{showlabels}
\allowdisplaybreaks

\newtheorem{mylemma}{Lemma}[section]
\newtheorem{mytheorem}{Theorem}[section]

\newtheorem{myremark}{Remark}[section]

\def\XXint#1#2#3{{\setbox0=\hbox{$#1{#2#3}{\int}$}
    \vcenter{\hbox{$#2#3$}}\kern-.5\wd0}}

\def\tx {\tilde{x}}

\def\tu {\tilde{u}}
\def\tphi {{\tilde{\phi}}}
\def\tT {\tilde{\cal T}}
\def\tcS {\tilde{\cal S}}
\def\tv {\tilde{v}}

\def\tP {\tilde{P}}
\def\tB {\tilde{B}}
\def\d{{\rm dist}}

\def\bI {{\bf I}}

\def\bA {{\bf A}}

\def\hG {\hat{G}}
\def\hu {\hat{u}}
\def\hw {\hat{w}}
\def\hf {\hat{f}}

\def\bi{{\bf i}}

\newcommand{\tcr}{\textcolor{red}}

\newcommand{\bv}{{\bf v}}
\newcommand{\bt}{{\bf t}}
\counterwithout{figure}{section}

\begin{document}
\title{Does PML exponentially absorb outgoing waves scattering from a periodic surface?}
\author{Wangtao Lu$^1$, Kuanrong Shen$^1$ and Ruming Zhang$^2$}
\footnotetext[1]{School of
  Mathematical Sciences, Zhejiang University, Hangzhou 310027, China. Email:
  wangtaolu@zju.edu.cn, kuanrongshen@zju.edu.cn. W.L. is partially supported by NSFC Grant
  12174310 and a Key Project of Joint Funds For Regional Innovation and
  Development (U21A20425).}
\footnotetext[2]{Institute of Mathematics, Technische Universit{\"a}t Berlin, Stra\ss{}e des 17. Juni 136, 10623 Berlin, Germany. Email: zhang@math.tu-berlin.de. The author is supported by the Deutsche Forschungsgemeinschaft (DFG, German Research Foundation) Project-ID 433126998.}
\maketitle
\begin{abstract}
The PML method is well-known for its exponential convergence rate and easy implementation for scattering problems with unbounded domains. For rough-surface scattering problems, authors in \cite{Chand2010} proved that the PML method converges at most algebraically in the physical domain. However, the authors also asked a question whether exponential convergence still holds for compact subsets.  In \cite{zha22}, one of our authors proved the exponential convergence for  $2\pi$-periodic surfaces via the Floquet-Bloch transform when $2k\in\mathbb{R}^+\backslash\mathbb{Z}$  where $k$ is the wavenumber; when $2k\in \mathbb{R}^+\cap\mathbb{Z}$, a nearly fourth-order convergence rate was shown in \cite{Zhang2022b}. The extension of this method to locally perturbed cases is not straightforward, since the domain is no longer periodic thus the Floquet-Bloch transform doesn't work, especially when the domain topology is changed.  Moreover, the exact decay rate when $2k\in \mathbb{R}^+\cap\mathbb{Z}$ remains unclear. The purpose of this paper is to address these two significant issues. For the first topic, the main idea is to reduce the problem by the DtN map on an artificial curve, then the convergence rate of the PML is obtained from the investigation of the DtN map. It shows exactly the same convergence rate as in the unperturbed case.   Second, to illustrate the convergence rate when $2k\in \mathbb{R}^+\cap\mathbb{Z}$, we design a specific periodic structure for which the PML converges at the fourth-order, showing that the algebraic convergence rate is sharp. We adopt a previously developed high-accuracy PML-BIE solver to exhibit this unexpected phenomenon.
\end{abstract}

\section{Introduction}

For the numerical simulation of wave scattering problems in unbounded domains, the perfectly matched layers (PML), which was invented by Berenger in  1994  
 in \cite{ber94}, is a widely used cut-off technique. The main idea of this method is to add an absorbing layer outside the physical domain, then the problem is approximated by a truncated problem with a proper { boundary} condition. We refer to \cite{joh08} to a detailed discussion of this method.
In this paper, we will focus on the PML convergence for  wave scattering problems with locally perturbed periodic surfaces.

In the past decades, many mathematicians have been working on the theoretical analysis and numerical implementations of scattering problems with (locally perturbed) periodic structures. For a special case, when the incident is quasi-periodic and the structure is purely periodic, there is a well established framework such that the problem is easily reduced into one periodicity cell. We refer to \cite{kir93,Bao1994a,Stryc1998} for theoretical discussions and \cite{Bao1995,Bao2005,Lechl2012} for numerical results.  However, when the surface is perturbed, or the incident field is not periodic (e.g., point sources, Herglotz wave functions), this framework no longer works. Despite the periodicity, the problem can be treated as a general rough surface scattering problem.  In \cite{Chand2005}, the well-posedness of the problems have been proved in a normal Sobolev space even for non-periodic surfaces and it is  also proved in weighted Sobolev spaces  in \cite{Chand2010}. For radiation conditions for the scattered field, we refer to \cite{hulurat21,hukir23}.

In order to take the advantage of the (locally perturbed) periodic structures, the Floquet-Bloch transform is applied to these problems, see \cite{Coatl2012,Hadda2016} for numerical implementations with absorbing media. In \cite{Lechl2015e}, the authors studied Herglotz wave functions scattered by periodic surfaces, with the help of the Floquet-Bloch transform. The method was extended to locally perturbed periodic surfaces in \cite{Lechl2016}. Based on these approaches, numerical methods have been developed. For 2D cases, the method was proved to be convergent in \cite{Lechl2016a,Lechl2017}. Based on  a detailed study of the singularity of the DtN map, a higher order method was also developed in \cite{Zhang2017e}. However, since the singularity of the DtN map becomes much more complicated for 3D cases, the convergence was not proved in general (see \cite{Lechl2016b}) and the extension of the method proposed in \cite{Zhang2017e} becomes impossible. In addition, the DtN map is a non-local boundary condition, which is difficult to be implemented numerically. Thus we are motivated to apply the PML method to solve these problems.

For the particular case that the incident field is quasi-periodic and the surface is purely periodic, exponential convergence for the PML method has been proved in \cite{chewu03} and the computation was carried out with an adaptive finite element method. For rough surface scattering problems, authors have proved in \cite{Chand2010} that the convergence in the whole physical domain is at most algebraic. At the end of this paper, the authors proposed a conjecture that exponential convergence also holds for compact subsets. In \cite{zha22}, { one of our authors proved the conjecture} for non-periodic incident fields scattered by purely periodic surfaces using the method of the Floquet-Bloch transform with a countable number of wavenumbers excluded; { at these excluded wavenumbers, the author further derived a nearly fourth-order decaying upper bound for the PML truncation error although its sharpness remains unjustified.} In this paper, we study the case with locally perturbed periodic surfaces. Since the domain is no longer periodic, the Floquet-Bloch transform no longer works. To this end, our idea is to first derive an upper bound for the difference of { two Dirichlet-to-Neumann (DtN) maps of both the original problem and the PML-truncated problems on a 
 common artificial curve separating an unperturbed region from the perturbed part. The DtN maps are constructed via single-layer operators defined through closely related transmission problems with the unperturbed periodic surface. The difference is estimated based on the theories in \cite{zha22,Zhang2022b}. The convergence result is then proved based on standard error estimates for elliptic differential equations. We construct a specific periodic structure to illustrate that the PML truncation error indeed can decay algebraically at a fourth order convergence rate. Finally, we adopt a recently developed high-accuracy PML-BIE method \cite{yuhulurat22} to numerically validate such an unexpected phenomenon.}

The rest of this paper is organized as follows. In section 2, the mathematical model for the problem is described and the main convergence results for the PML method are reviewed. In the next section, we study the convergence of the DtN map on an artificial curve. With this result, the local exponential convergence is proved in Section 4. We construct a fourth-order accurate PML in section 5. Numerical experiments are presented with the PML-BIE method in Section 6.

\section{Problem description}
\begin{figure}[htbp]
\centering
(a)\includegraphics[width=0.46\textwidth]{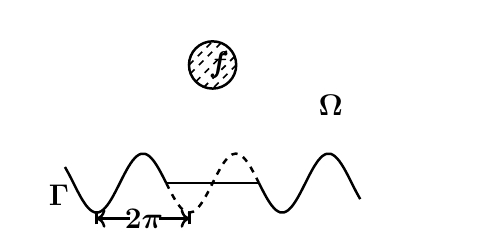}
(b)\includegraphics[width=0.46\textwidth]{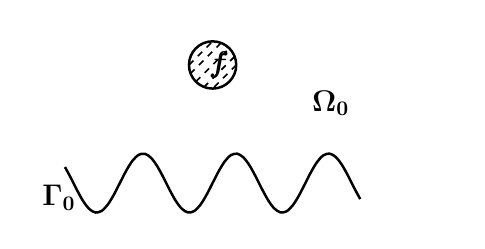}
\label{fig:model}
\caption{Profile of the scattering problem: $\Gamma$ differs from $\Gamma_0$ by the locally perturbed part; $f$ indicates the source term.}
\end{figure}

The profile of the scattering problem is depicted in Figure~\ref{fig:model}.
Let $\Gamma_0$ be a periodic and Lipschitz curve of period $2\pi$ in $x_1$-direction, where and in the following $(x_1,x_2)$ denotes the standard Cartesian coordinate system. Let a curve $\Gamma$ be from locally perturbing $\Gamma_0$ such that
$\Gamma\backslash\Gamma_0$ is bounded and Lipschitz (empty if
$\Gamma=\Gamma_0$). Let $\Omega$ and $\Omega_0$ be the two upper Lipschitz
domains bounded by $\Gamma$ and $\Gamma_0$, respectively. For technical reason,
we assume further that $\Omega$ (and hence $\Omega_0$) satisfies the following
geometrical condition
\begin{equation*}
  {\rm (GC):}\quad\quad (x_1,x_2)\in \Omega \implies (x_1,x_2+h) \in\Omega\quad \forall h\geq 0.
\end{equation*}

Consider the following scattering problem in the perturbed domain $\Omega$, 
\begin{align}
  \label{eq:gov:1}
  \Delta u + k^2 u =& f,\quad{\rm in}\quad\Omega,\\
  \label{eq:gov:2}
  u =& 0,\quad {\rm on}\quad\Gamma,
\end{align}
where $\Delta = \partial_{x_1}^2 + \partial_{x_2}^2$ denotes the 2D Laplacian,
$f\in \widetilde{H^{-1}}(\Omega)=[H^1(\Omega)]'$ is the exciting source term with a compact
support, $k>0$ denotes the wavenumber for the source, and $u$ denotes the
wavefield. Such a problem can describe a TE-polarized electric field with the
nonzero component $u$ propagating in $\Omega$ due to a perfectly electric
conductor $\Gamma$, or a sound field $u$ due to a sound-soft surface $\Gamma$.

To ensure the well-posedness of the above problem, one may enforce the following
upward propagating radiation condition (UPRC):
\begin{equation}
  \label{eq:uprc}
  u(x) = 2\int_{\Gamma_H} \frac{\partial G(x;y)}{\partial y_2} ds(y),
\end{equation}
where $x=(x_1,x_2)$, $y=(y_1,y_2)$, $\Gamma_H = \{(x_1,H):x_1\in\mathbb{R}\}$
denotes a straight line strictly above $\Gamma$ for some sufficiently large
$H>0$ such that ${\rm supp}f\subset \Omega_H:=\Omega\cap\{x: x_2<H\}$, and
  $G(x;y) = \frac{\bi}{4} H_0^{(1)}(k|x-y|)$
denotes the fundamental solution for
the Helmholtz equation (\ref{eq:gov:1}). Alternatively but more precisely, one
can enforce the half-plane Sommerfeld radiation condition (hpSRC):
  \begin{equation}
    \label{eq:hpsrc}
    \lim_{r\to\infty}\sup_{\alpha\in[0,\pi]}\sqrt{r}\left| \partial_{r}u(x)-\bi ku(x)\right| = 0,\ \sup_{r\geq R}r^{1/2}|u(x)|<\infty,\ {\rm and}\ u\in H_\rho^{1}(S_H^R),
  \end{equation}
  where $x=(r\cos\alpha,H+r\sin\alpha)$, $S_H^R=\{x\in\Omega:|x_1|>R, x_2<H\}$,
  and $H_\rho^1(\cdot)=(1+x_1^2)^{-\rho/2}H^1(\cdot)$ denotes a weighted Sobolev
  space. 

  \begin{myremark}
    In real applications, the total wave field $u$ is usually generated by
    specifying an incident plane or cylindrical wave in $\Omega$ and cannot be
    characterized directly by the above scattering problem. Nevertheless, as
    shown in \cite{hulurat21}, one can always decompose $u$ as the sum of a
    known (or more easily computed) function and another unknown wave field
    satisfying (\ref{eq:gov:1}), (\ref{eq:gov:2}), and one of the two radiation
    conditions (\ref{eq:uprc}) and (\ref{eq:hpsrc}). Thus, all results in this
    paper can be trivially extended to plane-wave or cylindrical-wave
    incidences.
  \end{myremark}

  We collect some well-known results in existing literature in the following.
  The well-posedness of the above scattering problem has been justified in
  \cite{hulurat21} as follows.
  \begin{mytheorem}
    \label{thm:wp}
    Under the geometrical condition (GC), the scattering problem
    (\ref{eq:gov:1}) and (\ref{eq:gov:2}) with either UPRC (\ref{eq:uprc}) or
    hpSRC (\ref{eq:hpsrc}) at infinity, has a unique solution $u\in
    H^{1}_{\rm loc}(\Omega)$ for any $f\in \widetilde{H^{-1}}(\Omega)$ with a compact support and any
    $k>0$. Moreover $u|_{\Omega_H}\in H_\Gamma^{1}(\Omega_H):=\{v\in
    H^1(\Omega_H): v|_{\Gamma}=0\}$.
  \end{mytheorem}

\begin{figure}[htbp]
\centering
 (a)\includegraphics[width=0.5\textwidth]{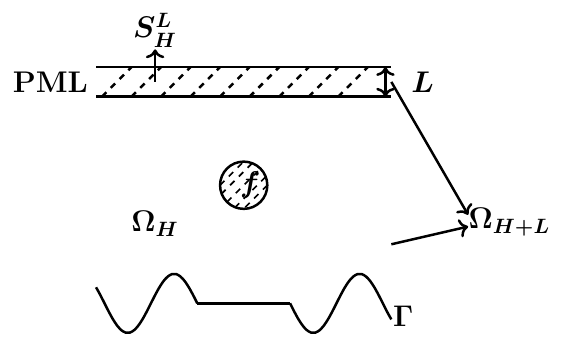}
(b)\includegraphics[width=0.4\textwidth]{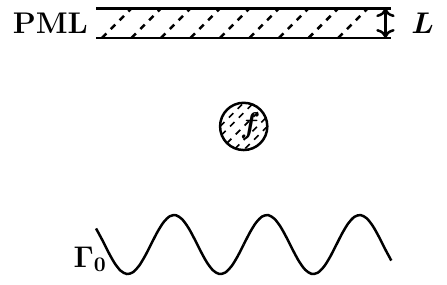}
\label{fig:pml}
\caption{Setup of the PML.}
\end{figure}

  To numerically solve the problem, it is advantageous to place a perfectly matched layer (PML)
  above $\Gamma_H$ to truncate $x_2$, as shown in Figure~\ref{fig:pml}.
  Mathematically, the PML can be characterized
  by a complexified coordinate transformation
  \begin{equation}
  \label{eq:sigma}
    \tx_2 = x_2 + \bi\int_0^{x_2}\sigma(t)dt,
  \end{equation}
      where $\sigma(x_2)=0$ for $|x_2|\leq H$ and $\sigma(x_2)>0$ for $|x_2|\in
  (H,\infty)$ that controls the absorption power of the PML. The planar strip $S_H^L=\{x: H< x_2 < H+L\}$ is called the PML
  region where $u(x)$ is analytically continued to $\tu(x)=u(x_1,\tx_2)$. On the
  PML boundary $\Gamma_{H+L}=\{x:x_2=H+L\}$, we assume that $\tu$, after going
  through the PML strip $S_H^L$, is absorbed sufficiently such that it is
  reasonable to assume $\tu=0$. Let $\Omega_{H+L}=\Omega\cap\{x:x_2<H+L\}$.
  Consequently, $\tu$ satisfies the following PML problem
  \begin{align}
    \label{eq:pmlgov:1}
    \nabla\cdot(\bA\nabla\tu) + k^2\alpha \tu =& f,\quad{\rm in}\quad\Omega_{H+L},\\
    \label{eq:pmlgov:2}
    \tu =& 0,\quad{\rm on}\quad \Gamma\cup\Gamma_{H+L},
  \end{align}
  where $\bA = {\rm diag}\{\alpha(x_2),\alpha^{-1}(x_2)\}$ and
  $\alpha(x_2)=1+\bi \sigma(x_2)$. Let $\tP = \int_H^{H+L}\alpha(x_2)d x_2$. The well-posedness of the
  PML problem has been justified in \cite{chamon09} as shown below.
  \begin{mytheorem}
    \label{thm:wp:pml}
    Under the geometrical condition (GC), the PML problem
    (\ref{eq:pmlgov:1}) and (\ref{eq:pmlgov:2}), has a unique solution $u\in
    H_0^{1}(\Omega_{H+L}):=\{v\in H^1(\Omega_{H+L}): v|_{\Gamma\cup\Gamma_{H+L}}=0\}$ for any $f\in \widetilde{H^{-1}}(\Omega)$ and any
    $k>0$ such that for any bounded domain $D\in\Omega_{H}$
    \[
    ||u||_{H^{1}(D)} \leq C ||f||_{\tilde{H}^{-1}(\Omega)}, 
    \]
    provided that $|\tP|$ is sufficiently large, where $C$ does not depend on $L$ or $\sigma$.
  \end{mytheorem}

  After the truncation of the vertical $x_2$-direction, $\Omega$ becomes a strip
  $\Omega_{H+L}$ which consists of two periodic semiwaveguides at infinity.
  Thus, existing techniques such as Floquet-Bloch transforms, Ricatti-equation
  governing marching operators, and recursive doubling procedures can be used to
  terminate the two semiwaveguides in terms of posing exact Neumann-to-Dirichlet
  or Dirichlet-to-Neumann maps. The resulting boundary value problem can then be
  solved by standard solvers. Naturally, one would ask how accurate would the
  PML-truncated solution $\tu$ be, compared with the total field $u$, in the
  physical domain $\Omega_{H}$.

  Chandler-Wilder and Monk \cite{chamon09} firstly studied the PML convergence
  theory when $\Gamma$ is a rough surface, not necessarily period at infinity.
  They proved that $\tu$ converges to $u$ at an algebraic rate in
  $\Omega_H$ but conjectured that an exponential rate in any compact subset of $\Omega_H$ as the
  PML parameter $|\tP|\to\infty$, and have strictly proved this for a
  flat $\Gamma$. In a recent work \cite{yuhulurat22}, Yu {\it et al.} proved an
  algebraically convergent rate of $\tu$ in $\Omega_H$ for the locally perturbed
  periodic curve $\Gamma$ under consideration and had numerically verified the
  exponentially convergent rate of $\tu$ on a compact subset of $\Omega_H$, yet
  a theoretical justification remains open.

  {For the unperturbed case $\Gamma=\Gamma_0$, One of our authors in
  \cite{zha22} adopted the Floquet-Bloch transform to firstly decouple the
  original problem into a series of subproblems, each of which possesses only
  quasi-periodic solutions, i.e., Bloch waves. The original problem, on the other hand, is written as the integral of these quasi-periodic solutions with respect to the quasi-periodicity parameters on a bounded interval. 
  Based on the contour deformation theory, the integral contour is modified near the Rayleigh anomalies and the locally exponential convergence
  of the PML-truncated solution $\tu$ to $u$ was successfully established for $2k\notin\mathbb{Z}$. If $2k\in\mathbb{Z}$, {the two Rayleigh anomalies coincide, making the contour deformation technique in \cite{zha22} break down}. Nevertheless, a high  order algebraically convergent can still be proved, with a detailed study of the convergence rate of the quasi-periodic PML solutions near the Rayleigh anomalies. We refer to \cite{Zhang2022b} for 3D bi-periodic cases and the method proposed is easily applied to 2D cases.  The two results are summarized below.
  \begin{mytheorem}
    \label{thm:main:unperturbed}
    If $\Gamma$ is purely periodic, i.e., $\Gamma=\Gamma_0$, and if the
    geometrical condition (GC) is satisfied, then: for any bounded open subset $D\subset \Omega_{H}$ and any $f\in \widetilde{H}^{-1}(\Omega)$,
    provided that $|\tP|$ is sufficiently large,
    \\
    (1) If $0<2k\notin\mathbb{Z}^+$, {there are constants $C,c_0>0$ such that}
    \begin{equation}
      ||\tu - u||_{H^1(D)}\leq C {e^{-c_0 |\tP|}}||f||_{\widetilde{H}^{-1}(\Omega)},
    \end{equation}
    (2) If $2k\in\mathbb{Z}^+$, {for any constant $0<\gamma_0<1$, there is a constant $C=C(\gamma_0)>0$ such that}
    \begin{equation}
      ||\tu - u||_{H^1(D)}\leq C  {|\tP|}^{-4\gamma_0}||f||_{\widetilde{H}^{-1}(\Omega)},
    \end{equation}
  \end{mytheorem}
  
The main contribution of this paper is to prove the exponential/algebraic convergence of
$\tu$ to $u$ in any compact subset of $\Omega_H$, when $\Gamma$ is a locally
pertubed periodic curve, i.e., $\Gamma\neq \Gamma_0$. 
{Note that the
  method of Floquet-Bloch transform in \cite{zha22} cannot be trivially extended to such a more
  general case since the Floquet-Bloch transform is not valid for non-periodic domains.} As the methodology does not rely on whether $k$ is a half integer or not, we shall assume $2k\notin\mathbb{Z}^+$ for the moment.
  
  The basic idea goes as follows. Firstly, for each of the original scattering problem and the PML problem, we shall, following the approach in \cite{yuhulurat22}, use a
single-layer operator to express the field in the exterior of a bounded domain
containing the perturbed part $\Gamma\backslash\Gamma_0$ so as to establish 
a Dirichlet-to-Neumann (DtN) map that truncates the unbounded problem into a
boundary value problem. {The single-layer operator is defined by an equivalent transmission problem involving only the unperturbed periodic surface.} Secondly, we shall justify
the exponentially decaying difference between the two DtN maps as $\tP$ increase. Finally, for the two boundary value problems, we present
the equivalent variational formulations with the aid of the two DtN maps, analyze
their inf-sup conditions, and establish the exponential convergence theory. 

\section{Dirichlet-to-Neumann maps}
As shown in Figure~\ref{fig:dtn},
\begin{figure}[htbp]
\centering
 (a)\includegraphics[width=0.8\textwidth]{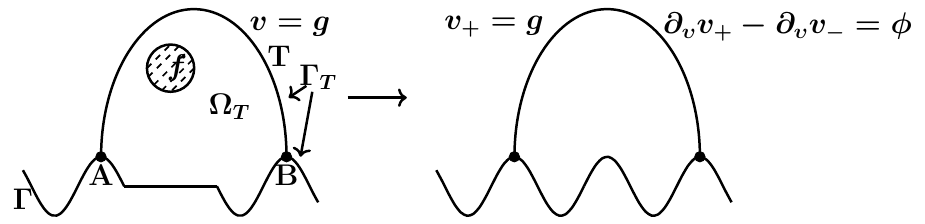}\\
(b)\includegraphics[width=0.8\textwidth]{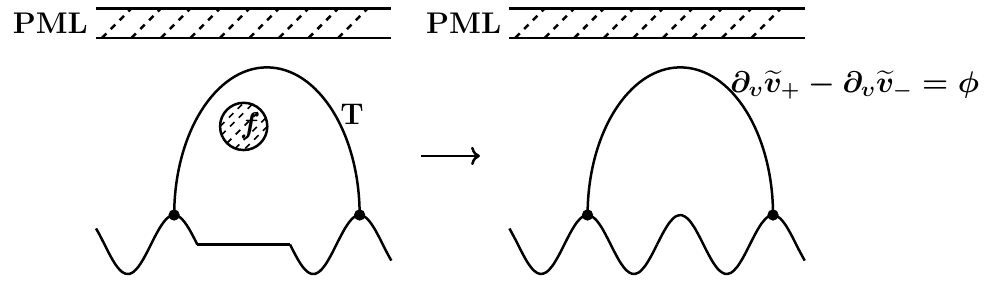}
\label{fig:dtn}
\caption{Dirichlet-to-Neumann maps for the original problem and the PML-truncating problem.}
\end{figure}
let $T$ be a sufficiently smooth curve between the
scattering surface $\Gamma$ and $\Gamma_H$, with endpoints $A$ and $B$ on
$\Gamma$. Let $\Omega_T$ be the domain bounded by $\Gamma$ and $T$ and let
$\Gamma_T = \overline{T}\cup(\Gamma\backslash\overline{\Omega_T})$ be the curve
consists of $T$ and the part of $\Gamma$ outside $\Omega_T$. We choose $T$
properly such that $\Gamma_T$ is Lipschitz and that $\Omega_T$ satisfies (GC).
For any given function $g\in \widetilde{H^{1/2}} (T)=[H^{-1/2}(T)]'\subset H^{1/2}(\Gamma_T)$,
consider the following two problems: find $v\in H^{1}_{\rm
  loc}(\Omega\backslash\overline{\Omega_T})$ such that
\begin{align*}
  (P1):\quad\quad \left\{
  \begin{array}{ll}
  \Delta v + k^2 v = 0,&\quad{\rm in}\quad \Omega\backslash\overline{\Omega_T},\\
  v = g,&\quad{\rm on}\quad \Gamma_T,\\
  v\ {\rm satisfies\ (\ref{eq:uprc})\ or\ (\ref{eq:hpsrc}),}
  \end{array}
       \right.
\end{align*}
and $\tv\in H^1(\Omega_{H+L}\backslash\overline{\Omega_T})$ such that
\begin{align*}
  (P2):\quad\quad\left\{
  \begin{array}{ll}
  \nabla\cdot(\bA\nabla\tv)  + k^2\alpha \tv = 0,&\quad{\rm in}\quad \Omega_{H+L}\backslash\overline{\Omega_T},\\
  \tv = g,&\quad{\rm on}\quad \Gamma_T,\\
  \tv = 0,&\quad{\rm on}\quad \Gamma_{H+L}.
  \end{array}
 \right.
\end{align*}
One can find a function $f_g\in \widetilde{H^{1}}(\Omega_{H})$ with its compactly support in $\Omega_H$ such that $f_g|_{T}=g$. Then, the two functions $v-f_g$ and $\tv-f_g$ respectively satisfy the original scattering problem and the PML problem with $f$ replaced by $\Delta f_g + k^2 f_g\in
\widetilde{H^{-1}}(\Omega_{H})$. According to Theorems~\ref{thm:wp} and
\ref{thm:wp:pml}, both (P1) and (P2) are well-posed. Thus, we can define two DtN
maps ${\cal T}, \tT: \widetilde{H^{1/2}}(T)\to H^{-1/2}(T)$ that are bounded  and satisfy ${\cal T}v|_{T}=\partial_{\nu}v|_{T}$ and $\tT \tv|_T =
 \partial_{\nu}{\tv}|_{T}$ where $\nu$ denotes the outer unit normal to $T$. Clearly,
the two DtN maps ${\cal T}$ and $\tT$ can respectively serve as the transparent
boundary conditions of the scattering problem and the PML problem to
truncate them into boundary value problems. The purpose of this section
is to estimate the difference between ${\cal T}$ and $\tT$ when $\tP$ is sufficiently large. To study this, we adopt the
idea in \cite{yuhulurat22} and \cite{hukir23} to use single-layer operators to
express the two DtN maps.

Consider two associate problems in the unperturbed domain $\Omega_0$ and the
PML region $\Omega_{0,H+L}:=\Omega_0\cap \{x:x_2<H+L\}$: Given
$\phi,\tphi\in H^{-1/2}(T)$, find $v\in H^{1}_{\rm loc}(\Omega_0)$ such that
\begin{align*}
  (P1'):\quad\quad \left\{
  \begin{array}{ll}
  \Delta v + k^2 v = 0,&\quad{\rm in}\quad \Omega_0\backslash{\overline{T}},\\
    \partial_{\nu }v_+ - \partial_{\nu }v_- = \phi, &\quad{\rm on}\quad T,\\
  v = 0,&\quad{\rm on}\quad \Gamma_0,\\
  v\ {\rm satisfies\ (\ref{eq:uprc})\ or\ (\ref{eq:hpsrc}),}
  \end{array}
       \right.
\end{align*}
and $\tv\in H^1(\Omega_{0,H+L})$ such that
\begin{align*}
  (P2'):\quad\quad\left\{
  \begin{array}{ll}
  \nabla\cdot(\bA\nabla\tv)  + k^2 \alpha\tv = 0,&\quad{\rm in}\quad \Omega_{0,H+L}\backslash{\overline{T}},\\
    \partial_{\nu }\tv_+ - \partial_{\nu }\tv_- = \tphi, &\quad{\rm on}\quad T,\\
  \tv = 0,&\quad{\rm on}\quad \Gamma_{H+L}\cup \Gamma_0,
  \end{array}
 \right.
\end{align*}
where the two subscripts $+$ and $-$ indicate that the normal derivatives are
taken from the exterior and interior of $\Omega_T$, respectively. 

Let $\Omega_{0,H}=\Omega_0\cap\{x:x_2<H\}$. Then, one can construct functions
$f_\phi,f_\tphi\in \widetilde{H^1}(\Omega_{0,H})$ such that $\partial_{\nu
}[f_\phi]_+ - \partial_{\nu }[f_\phi]_- = \phi$ and $\partial_{\nu
}[f_{\tphi}]_+ - \partial_{\nu }[f_{\tphi}]_- = \tphi$. We also assume that there is a constact $C>0$ such that $\|f_\phi\|_{\widetilde{H^1}(\Omega_{0,H})}\leq C\|\phi\|_{H^{-1/2}(T)}$ and $\|f_\tphi\|_{\widetilde{H^1}(\Omega_{0,H})}\leq C\|\tphi\|_{H^{-1/2}(T)}$.
Again, $v - f_\phi$ and
$\tv - f_{\tphi}$ respectively satisfy the original scattering problem and the
PML problem with $f$ replaced by the two functions $-\Delta
f_\phi-k^2f_\phi$ and {$-\Delta f_\tphi-k^2f_\tphi$ in $
\widetilde{H^{-1}}(\Omega_H)$}.  Thus, $(P1')$ and $(P2')$ (for sufficiently large
$L$ and $\tP$) are well-posed. Consequently, we can define two bounded 
operators ${\cal S}, \tcS: H^{-1/2}(T) \to \widetilde{H^{1/2}}(T)$ such that
${\cal S}\phi = v|_{T}$ and $\tcS\tphi = \tv|_{T}$. 
{Using the background Green functions for the two problems, it is easy to see that the two operators coincide with the standard single-layer operators.}
Before proceeding, we present some important properties of the single-layer operator ${\cal S}$ defined on $T$.
\begin{mylemma}
 For any $k>0$, the operator ${\cal S}$ is Fredholm of index zero. Moreover, there exists a smooth curve $T$ such that ${\rm supp}f \subset \Omega_T$, $\Gamma_T$ satisfies (GC), and that ${\cal S}$ is boundedly invertible.
 \begin{proof}
 The proof follows simply from the two works \cite{hulurat21,hukir23}. Consider problem 
 $(P1')$ for $k=\bi$. The variational solution $\tv_\bi$ of $(P1')$ for $k=\bi$ belongs to
{$H^1_{(\tau)}(\Omega_{0})\!:=\!\{u\!\in\! H_0^1(\Omega_{0})\!: e^{\tau |x|} u\!\in\!
H^1(\Omega_{{0}})\}$ for $\tau\in(0,1)$,
and the mapping $\phi\,\mapsto \tv_\bi$ is bounded from $H^{-1/2}(T)$
into $H^1_{(\tau)}(\Omega_{0})$ and even compact from
$H^{-1/2}(T)$ into $L^2_{(\tau')}(\Omega_{{0}}):=\{u\in
H_0^1(\Omega_{{0}}): e^{\tau' |x|} u\in L^2(\Omega_{{0}})\}$
 with $\tau'<\tau$ (cf.~\cite[Theorem 3.1]{hukir23}).
 Denote the operator mapping $\phi\,\mapsto
 \tilde{v}_\bi|_{T}$ by ${\cal S}_{\bi}$}. The invertibility of ${\cal S}_\bi$ follows from the strong ellipticity of $(\Delta-1)$. 
 
It follows that
 $({\cal S}-{\cal S}_{\bi})\phi=w|_{T}$ where $w\in H_{0,{\rm loc}}^1(\Omega_0)$ is
 the solution to
\begin{equation}\label{eq:w}
\Delta w + k^2 w = {-(k^2+1)} \tilde{v}_\bi,
\quad{\rm in}\quad{\Omega}_{{0}}.
\end{equation}
By \cite[Theorem 3.3]{hukir23}, the unique solvability of $w$ to the previous equation \eqref{eq:w} yields the
boundedness of the mapping 
$L^2_{(\tau')}(\Omega_{0})\ni \tilde{v}_\bi \mapsto w\in
\widetilde{H^{1/2}}(T)$.  
This together with the compactness of
$\phi\,\mapsto \tilde{v}_\bi$ from $H^{-1/2}(T)\rightarrow
L^2_{(\tau')}(\Omega_{0})$ proves the compactness of ${\cal
  S}-{\cal S}_{\bi}$.
 
 To prove the invertibility of ${\cal S}$, we justify that equation ${\cal S}\phi = v|_{T} = 0$ possesses the zero solution only. Let $\Omega_-$ be the domain bounded by $T$ and $\Gamma_0$ and $\Omega_+=\Omega_0\backslash\overline{\Omega_-}$. Since $v|_{T}= 0$, $v|_{\Omega_+}$ solves the original scattering problem (\ref{eq:gov:1}) and (\ref{eq:gov:2}) with $f\equiv 0$ and $\Gamma$ replaced by $\partial\Omega_+$. Thus, $v|_{\Omega_+}\equiv 0$ so that $v|_{\Omega_-}$ solves 
 \begin{align*}
     -\Delta v =& k^2 v,\quad {\rm in}\quad\Omega_-,\\
     v=&0,\quad{\rm on}\quad \partial\Omega_-.
 \end{align*}
 If $k^2$ is not a Dirichlet eigenvalue of $-\Delta$ on $\Omega_-$, then $v|_{\Omega^-}=0$ so that 
 \[
 \phi = \partial_{\nu}v_+ - \partial_{\nu} v_-=0,
 \]
 which justifies the injectivity and hence the bijectivity of ${\cal S}$.
 
Suppose now $k^2$ is a Dirichlet eigenvalue of $-\Delta$ for the specified curve $T$. Choose another curve $T_0$ intersecting $\Gamma_0$ at the endpoints of $T$ such that the vertical distance between $T_0$ and $\Gamma_0$ is sufficiently small. Then, Fredrich's inequality implies that $k^2$ cannot be a Dirichlet eigenvalue of $-\Delta$
for the domain bounded by $\Gamma_0$ and $T_0$. Let $\xi\in[0,1]$ and consider the following family of curves
\[
T(\xi):=\{(x_1,\xi x_2^0(x_1)+(1-\xi)x_2^1(x_1)):(x_1,x_2^0(x_1))\in T_0, (x_1,x_2^1(x_1))\in T\}.
\]
The corresponding sesquilinear form of $\Delta u+k^2 u$ for the domain $\Omega(\xi)$ bounded by $\Gamma_0$ and $T(\xi)$ defines a linear operator ${\cal L}(\xi): H_0^{1}(\Omega(\xi))\to [H_0^{1}(\Omega(\xi))]^*$ that analytically depends on $\xi$ and is Fredholm of index zero. Since ${\cal L}(1)$ is invertible, there exists at most countable values of $\xi$ in $[0,1]$ such that ${\cal L}(\xi)$ has a non-degenerate kernel. Therefore, for a sufficiently small parameter $\epsilon>0$, there exist $\xi_0\in(0,\epsilon)$  such that ${\rm supp}f\in \Omega(\xi_0)$ and that ${\cal L}(\xi_0)$ is invertible, i.e., $k^2$ is not an eigenvalue of $-\Delta$ when $T$ is deformed to $T(\xi_0)$. 
 \end{proof}
\end{mylemma}

Choosing $\tphi=\phi$, we have $f_\phi=f_\tphi$ so that when $2k\notin\mathbb{Z}$,
  Theorem~\ref{thm:main:unperturbed} implies 
  \begin{align}
  \label{eq:diff:tvv}
  ||\tv - v||_{H^1(D)}&\leq C e^{-c_0  |\tP|}||-\Delta f_\phi-k^2f_\phi||_{\widetilde{H}^{-1}(\Omega_H)}\nonumber\\
  &\leq Ce^{-c_0 |\tP|} ||\phi||_{H^{-1/2}(T)},
\end{align}
for any bounded domain $D\subset \Omega_0$. Choosing $D$ sufficiently large to
contain $T$, we have by the trace theorem that 
\begin{align*}
  ||(\tcS-{\cal S})\phi|| =&||(\tv-v)|_T|| \leq C||\tv - v||_{H^1(D)}\\
  \leq& Ce^{-c_0 |\tP|} ||\phi||_{H^{-1/2}(T)},
\end{align*}
implying $||\tcS - {\cal S}||\leq Ce^{-c_0 |\tP|}$.
They indicate the following lemma.
\begin{mylemma}
  \label{lem:conv:tcS}
  Suppose $2k\notin\mathbb{Z}.$ For sufficiently large $|\tP|$, $\tcS$ is boundedly invertible and   
  \begin{equation}
    ||\tcS^{-1} - {\cal S}^{-1}|| \leq  C e^{-c_0 |\tP|}.
  \end{equation}

  \begin{proof}
First, choose sufficiently large $|\tP|$ such that
\begin{align*}
  ||\tcS {\cal S}^{-1} - {\cal I}||\leq  Ce^{-c_0 |\tP|} ||{\cal S}^{-1}|| < 1/2.
\end{align*}
 By the method of Neumann series, $\tcS{\cal S}^{-1}$ as a map from
 $H^{1/2}(T)$ to itself  is invertible and   
 \begin{align*}
   ||(\tcS {\cal S}^{-1})^{-1}||\leq 2.
 \end{align*}
 One directly verifies that $\tcS^{-1} = {\cal S}^{-1}(\tcS {\cal
   S}^{-1})^{-1}$ is bounded invertible with $||\tcS^{-1}||\leq 2 ||{\cal S}^{-1}||$ and that 
 \begin{align*}
||\tcS^{-1}-{\cal S}^{-1}|| = ||\tcS^{-1}({\cal I}-\tcS{\cal S}^{-1})||\leq 2 C||{\cal S}^{-1}||^2 e^{-c_0|\tP|}.
 \end{align*}
  \end{proof}
\end{mylemma}

Now, for any $g\in \widetilde{H^{1/2}}(T)$, set $\phi = {\cal S}^{-1}g$ in
$(P1')$ and $\tphi = \tcS^{-1}g$ in $(P2')$. 
Then, ${\cal T} g = \partial_{\nu} v_+$ and $\tilde{\cal T} g = \partial_{\nu} \tilde{v}_+$.
Lemma~\ref{lem:conv:tcS}
implies 
\begin{align}
  ||\tphi-\phi||_{H^{-1/2}(T)} \leq Ce^{-c_0 |\tP|} ||g||_{\in \widetilde{H^{1/2}}(T)}.
\end{align}

Clearly, the restriction of the solution $v$ of $(P1')$ onto
$\Omega\backslash\overline{\Omega_T}$ is the solution of $(P1)$, and the
restriction $\tv|_{\Omega_{H+L}\backslash\overline{\Omega_T}}$ the solution of
$(P2)$. We compare $\partial_\nu v_+$ and $\partial_{\nu} \tv_+$ on $T$.
Decomposing $\phi = \tphi + (\phi-\tphi) = \tcS^{-1} g + ({\cal
  S}^{-1}-\tcS^{-1})g$ in $(P1')$, we see that $v$ is the sum of $v^1$ where
$\phi$ is replaced by $\tcS^{-1}g$ in $(P1')$ and $v^2$ where $\phi$ is replaced by $({\cal S}^{-1}-\tcS^{-1})g$. Moreover,
choosing the bounded domain $D$ large enough to contain $T$, by (\ref{eq:diff:tvv}),
\begin{align*}
  ||\partial_\nu v^1_+ - \partial_\nu \tv_+||_{H^{-1/2}(T)}\leq C || v^1 - \tv||_{H^1(D)} \leq Ce^{-c_0 |\tP|} ||g||_{\widetilde{H^{1/2}}(T)}.
\end{align*}
By Lemma~\ref{lem:conv:tcS} and by the well-posedness of $(P1')$,
\begin{align*}
  ||\partial_\nu v^2_+||_{H^{-1/2}(T)}\leq C || v^2||_{H^1(D)} \leq  C || ({\cal
  S}^{-1}-\tcS^{-1})g||_{H^{-1/2}(T)}\leq Ce^{-c_0 |\tP|} ||g||_{\widetilde{H^{1/2}}(T)}.
\end{align*}
The triangular inequality then implies
\begin{align*}
  ||({\cal T} - \tT) g||_{H^{-1/2}(T)}=||\partial_\nu v_+ - \partial_\nu \tv_+||_{H^{-1/2}(T)}\leq Ce^{-c_0 |\tP|} ||g||_{\widetilde{H^{1/2}}(T)},
\end{align*}
giving rise to 
\begin{mylemma}
\label{lem:diffTtT}
  Suppose $2k\notin\mathbb{Z}.$ For sufficiently large $\tP$, 
  \begin{equation}
    ||{\cal T} - \tT|| \leq  C e^{-c_0 |\tP|}.
  \end{equation}
\end{mylemma}

\section{Local convergence of the PML solution}
Let $\Gamma_p = \partial\Omega_T\cap\Gamma\backslash\overline{T}$ be the part of $\Gamma$ bounded by the two endpoints of $T$. With the two DtN maps ${\cal T}$ and $\tT$ well-defined, the original scattering problem (\ref{eq:gov:1}), (\ref{eq:gov:2}) and (\ref{eq:hpsrc}) can now be truncated as the following boundary value problem
\begin{equation*}
    ({\rm OP}): \left\{
    \begin{array}{ll}
    \Delta u + k^2 u = f,&{\rm in}\quad \Omega_T, \\
    u = 0,&{\rm on}\quad \Gamma_p, \\
    \partial_{\nu} u = {\cal T} u,&{\rm on}\quad  T,
    \end{array}
    \right.
\end{equation*}
and the PML problem (\ref{eq:pmlgov:1}) and (\ref{eq:pmlgov:2}) as 
\begin{equation*}
    ({\rm TP}): \left\{
    \begin{array}{ll}
    \Delta \tu + k^2 \tu = f,&{\rm in}\quad \Omega_T, \\
    \tu = 0,&{\rm on}\quad \Gamma_p, \\
    \partial_{\nu} \tu = \tT \tu,&{\rm on}\quad  T. \\
    \end{array}
    \right.
\end{equation*}
We note that in Problem (TP), equation (\ref{eq:pmlgov:1}) reduces to the original Helmholtz equation since the computational region $\Omega_T$ is away from the PML region such that ${\bA} = \bI$ and $\alpha = 1$. Now, we consider the variational formulations of the two problems. Let $V:=\{v\in H^1(\Omega_T): v|_{\Gamma_T\cap \partial\Omega_T} = 0\}$, and $a,\tilde{a}: V\times V\to \mathbb{C}$ be two bilinear forms given by 
\begin{align}
    a(p,q) =& -(\nabla p,\nabla q)_{\Omega_T} + k^2(p,q)_{\Omega_T} + \langle {\cal T} \gamma_T p, \gamma_T q\rangle_{T},  \\
    \tilde{a}(p,q) =& -(\nabla p,\nabla q)_{\Omega_T} + k^2(p,q)_{\Omega_T} + \langle \tT \gamma_T p, \gamma_T q\rangle_{T},  
\end{align}
where $(\cdot,\cdot)_{\Omega_T}$ denotes the standard $L^2$ inner product, $\langle \cdot,\cdot\rangle_{T}$ denotes the duality pair between $H^{-1/2}(T)$ and $\widetilde{H^{1/2}}(T)$, and $\gamma_T: V\to \widetilde{H^{1/2}}(T)$ denotes the trace operator.
Clearly, (OP) is equivalent to the following weak formulation: Find $u\in V$, such that  
\begin{equation*}
    \noindent{\rm (WOP):}\quad\quad\quad a(u,v) := \langle f,v\rangle_{\Omega_T},\quad \forall v\in V,
\end{equation*}
where $\langle\cdot,\cdot\rangle_{\Omega_T}$ denotes the duality pair between $\widetilde{H^{-1}}(\Omega_T)$ and $H^1(\Omega_T)$.
Similarly, (TP) is equivalent to the following weak formulation: Find $\tu\in V$, such that  
\begin{equation*}
    \noindent{\rm (WTP):}\quad\quad\quad\tilde{a}(\tu,v) := \langle f,v\rangle_{\Omega_T},\quad \forall v\in V.
\end{equation*} 
The two sesqui-linear forms $a$ and $\tilde{a}$ induce two bounded operators $B, \tB: V\to V^*$ such that
\[
\left<Bp,q\right>_{\Omega_T} = a(p,q),\quad \left<\tB p,q\right>_{\Omega_T} = \tilde{a}(p,q),
\]
for all $p,q\in V$. Authors in \cite{hukir23} proved that $B$ in fact is Fredholm. Thus, 
the uniqueness of problem (WOP) implies that $B$ is bijective 
so that $||B^{-1}||\leq c^{-1}$ for some positive constant $c$. Let $\delta B= B - \tB$. By Lemma~\ref{lem:diffTtT}, for $2k\notin \mathbb{Z}$ and for a sufficiently large $L$ or $\tP$, 
\begin{align*}
|\left<\delta B p, q\right>_{\Omega_T}|  =& |a(p,q) - \tilde{a}(p,q)|\\
=&\left|\langle({\cal T} - \tT)\gamma_Tp,\gamma_T q\rangle\right|\\
\leq & C e^{-c_0 |\tP|} ||\gamma_Tp||_{\widetilde{H^{1/2}}(\Omega_T)}||\gamma_Tq||_{\widetilde{H^{1/2}}(\Omega_T)}\\
\leq & C e^{-c_0 |\tP|} ||p||_{H^1(\Omega_T)}||q||_{H^1(\Omega_T)},
\end{align*}
so that $||\delta B||\leq C e^{-c_0 |\tP|}$. Thus, the method of Neumann series indicates that $\tB$ is also bounded invertible as long as $||\delta B||< c$. Similar to the proof of Lemma~\ref{lem:conv:tcS}, we derive that 
\[
||\tB^{-1} - B^{-1}|| \leq Ce^{-c_0 |\tP|}.
\]
Consequently, Problem (WTP) has a unique solution $\tu\in V$ and 
\[
||u - \tu||_{H^1(\Omega_T)} = ||B^{-1} f - \tB^{-1} f||_{H^1(\Omega_T)} \leq  C e^{-c_0 |\tP|} ||f||_{\widetilde{H^{-1}(\Omega)}}.
\]
The above indicates that the PML solution $\tu$ converges to $u$ exponentially in the compact domain $\Omega_T$. 

We now claim that such a locally exponential convergence holds for any compact subdomain of the physical domain $\Omega_H$. 
Set $\phi = {\cal S}^{-1}\gamma_T u$ in
$(P1')$ and $\tphi = \tcS^{-1}\gamma_T \tu$ in $(P2')$. 
Clearly,  the restriction of the solution $v$ of $(P1')$ onto $\Omega\backslash\overline{\Omega_T}$ is $u|_{\Omega\backslash\overline{\Omega_T}}$, and $\tv|_{\Omega_{H+L}\backslash\overline{\Omega_T}}=\tu|_{\Omega_{H+L}\backslash\overline{\Omega_T}}$. Moreover, 
\begin{align*} ||\phi - \tphi||_{H^{-1/2}(T)}  \leq & ||({\cal S}^{-1}-\tcS^{-1})\gamma_T u ||_{H^{-1/2}(T)}+|| {\tcS}^{-1}\gamma_T (u - \tu) ||_{H^{-1/2}(T)}\\ \leq& C e^{-c_0 |\tP|} ||f||_{\widetilde{H^{-1}(\Omega)}}.
\end{align*}
As in section 3, 
we find two functions $f_\phi$ and $f_\tphi$  in $\widetilde{H^1}(\Omega_{0,H})$ such that $\partial_{\nu
}[f_\phi]_+ - \partial_{\nu }[f_\phi]_- = \phi$, $\partial_{\nu
}[f_{\tphi}]_+ - \partial_{\nu }[f_{\tphi}]_- = \tphi$, and $||f_\phi - f_\tphi||\leq C||\phi - \tphi||$. Applying Theorem~\ref{thm:wp:pml} (with $f$ replaced by $(-\Delta - k^2)[f_\phi-f_\tphi]$) and then Theorem~\ref{thm:main:unperturbed} (with $f$ replaced by $(-\Delta-k^2) f_\phi\in \widetilde{H^{-1}}(\Omega_H)$), 
we obtain the exponential convergence of the PML solution for any bounded domain exterior of $\overline{\Omega_T}$. Consequently, we obtain
\begin{mytheorem}
    \label{thm:main}
    {Suppose $2k\in\mathbb{R}^+\backslash\mathbb{Z}$}. Under the
    geometrical condition (GC) and provided that $\tP$ is sufficiently
    large,
    \begin{equation}
      ||\tu - u||_{H^1(D)}\leq C e^{-c_0 |\tP|}||f||_{\widetilde{H}^{-1}(\Omega)} 
    \end{equation}
    for any bounded open subset $D\subset \Omega_{H}$.
\end{mytheorem}

Following exactly the same procedure above, we obtain the algebraic convergence when $k$ is a half integer as stated below.
 \begin{mytheorem}
    \label{thm:main2}
    {Suppose $2k\in\mathbb{Z}^+$}. Under the
    geometrical condition (GC) and provided that $|\tP|$ is sufficiently
    large, for any fixed $\gamma_0\in(0,1)$,
    \begin{equation}
      ||\tu - u||_{H^1(D)}\leq C |\tP|^{-4\gamma_0}||f||_{\widetilde{H}^{-1}(\Omega)} 
    \end{equation}
    for any bounded open subset $D\subset \Omega_{H}$.
\end{mytheorem}

\begin{myremark}
   In \cite{zha22,Zhang2022b}, one of our authors directly applied the Floquet-Bloch transform 
    to establish the same PML convergence theory for purely periodic surfaces. The method is extendable to bounded penetrable medium or locally perturbed periodic surfaces, following the domain transformation method proposed in \cite{Lechl2017}.     Nevertheless, when the perturbation changes the topology, for example when the scattering domain contains an impenetrable obstacle, the method is no longer valid since the domain transformation can no longer be constructed. In comparison, our method is still extendable to all of the three situations provided that the original scattering problem is well-posed, and the rest is just a routine work as proposed in this paper.
\end{myremark}

\section{A fourth-order convergent PML}
We now study a specific example to illustrate that PML absorbs outgoing waves at most fourth order such that the convergence rate in Theorem~\ref{thm:main2}, is nearly sharp. Let us consider the following problem:
\begin{align}\label{eq:4th1}
\Delta u + (k^2+\epsilon \sin x_1 h(x_2))  u =& -\delta(x - x^*),\quad {\rm in}\quad\mathbb{R}_+^2,\\\label{eq:4th2}
u =& 0,\quad{\rm on}\quad \Gamma,
\end{align}
where $0<\epsilon\ll 1$,
\[
h(x_2) = \begin{cases}
    1 & x_2\in(0,1);\\
    0 & {\rm otherwise},
\end{cases}
\]
and $x^* = (x_1^*,x_2^*)$ with $x_2^*>1$. Without loss of generality, we assume $x_1^*=0$. 
\begin{myremark}
Instead of studying a scattering problem with a periodic surface, we here consider a periodic layered structure in the half space $\mathbb{R}^2_+$. This approach is to simplify the representation by avoiding the huge computational complexity brought by the domains transformations from the periodic domain to $\mathbb{R}^2_+$. However, the idea is extended without any difficulty for  periodic surface scattering problems.
\end{myremark}

From the perturbation theory, the well-posedness of the problem \eqref{eq:4th1}-\eqref{eq:4th2} is ensured for sufficiently small $\epsilon$. Moreover, it is easy to see that the solution is analytic w.r.t. $\epsilon$ at $\epsilon = 0$.  Thus, we suppose 
\[
u = \sum_{j=0}^{\infty} u_j\epsilon^j.
\]
The leading term $u_0$ is Green's function of the half-space $\mathbb{R}^2_+$ given by
\begin{align}
\label{eq:u0}
 u_0(x) =& \frac{\bi}{4} \left[H_0^{(1)}(k|x - x^*|) - H_0^{(1)}(k|x + x^*|)\right]\nonumber\\
 =& \frac{1}{4\pi}\int_{-\infty}^{+\infty}\frac{e^{-\bi\xi x_1}[e^{\bi \mu(\xi) |x_2-x_2^*|}-e^{\bi \mu(\xi)(x_2+x_2^*)}]}{\mu(\xi)}d\xi,
\end{align}
where $\mu(\xi) = \sqrt{k^2-\xi^2}$ with the negative real axis as its branch cut.
The second term $u_1$ is governed by the following source problem:
\begin{align}
\label{eq:gov:u1}
    \Delta u_1 + k^2 u_1 =& -\sin x_1 h(x_2) u_0,\quad x_2 > 0\\
\label{eq:bc:u1}
    u_1 =& 0,\quad {\rm on}\quad x_2 = 0.
\end{align}
In the following, let $\hat{f}$ denote the Fourier transform of a generic function $f(x_1)$ w.r.t. $x_1$ given by
\[
\hf(\xi) = \int_{-\infty}^{+\infty} f(x_1) e^{\bi \xi x_1} dx_1.
\]
Taking the $x_1$-Fourier transform of equations \eqref{eq:gov:u1} and \eqref{eq:bc:u1}, 
\begin{align}
    \hu''_1(x_2;\xi) + \mu^2(\xi) \hu_1(x_2;\xi) =& 
    \frac{\bi }{2}h(x_2) \left[\hu_0(x_2;\xi+1) - \hu_0(x_2;\xi-1)\right],\quad x_2 > 0,\\
    \hu_1(x_2;\xi) =& 0,\quad {\rm on}\quad x_2 = 0.
\end{align}
By \eqref{eq:u0}, 
\[
\hu_0(x_2;\xi) = \frac{1}{2}\frac{e^{\bi \mu(\xi) |x_2-x_2^*|}-e^{\bi \mu(\xi)(x_2+x_2^*)}}{\mu(\xi)}.
\]
For $x_2>1$ such that $h(x_2)\equiv 0$, we assume
\[
\hu_1(x_2;\xi)  = A(\xi) e^{\bi \mu(\xi) (x_2 - 1)},\quad x_2 > 1.
\]
For $x_2<1$ such that $h(x_2)\equiv 1$ and $x_2^*>x_2$, \eqref{eq:u0} implies 
\[
\hu_0(x_2;\xi) = \frac{1}{2}\frac{e^{\bi \mu(\xi) |x_2-x_2^*|}-e^{\bi \mu(\xi)(x_2+x_2^*)}}{\mu(\xi)} = -\bi e^{\bi \mu(\xi)x_2^*} \frac{\sin(\mu(\xi) x_2)}{\mu(\xi)}.
\]
Thus,
\[
\hu_1^-(x_2;\xi) = \frac{1}{2}\left[\frac{\sin(\mu(\xi+1)x_2)}{\mu(\xi+1)(2\xi+1)} e^{\bi \mu(\xi+1) x_2^*}   + \frac{\sin(\mu(\xi-1)x_2)}{\mu(\xi-1)(2\xi-1)} e^{\bi \mu(\xi-1) x_2^*}\right]
\]
provides a special solution such that we assume
\[
\hu_1(x_2;\xi) = \hu_1^-(x_2;\xi) + B(\xi) \sin(\mu(\xi) x_2), x_2 < 1.
\]
The continuity condition on $x_2=1$ leads to the following linear system
\begin{align*}
    A(\xi) =& \hu_1^-(1;\xi) + B(\xi) \sin(\mu(\xi)),\\
    \bi A(\xi)  =& \frac{[\hu_1^-]'(1;\xi)}{\mu(\xi)} + B(\xi) \cos(\mu(\xi)).
\end{align*}
On solving the linear system, we obtain for $x_2>1$ that
\begin{align}
  \label{eq:hu1}
\hu_1(x_2;\xi) 
=&\frac{1}{2}\Bigg\{
\frac{e^{\bi\mu(\xi+1) x_2^*}}{2\xi+1}\left[
\frac{\sin(\mu(\xi+1))}{\mu(\xi+1)}\cos(\mu(\xi)) -\frac{\sin(\mu(\xi))}{\mu(\xi)}\cos(\mu(\xi+1))\right]\nonumber\\
&+\frac{e^{\bi\mu(\xi-1) x_2^*}}{2\xi-1}\left[
\frac{\sin(\mu(\xi-1))}{\mu(\xi-1)}\cos(\mu(\xi)) -\frac{\sin(\mu(\xi))}{\mu(\xi)}\cos(\mu(\xi-1))\right]
\Bigg\}e^{\bi \mu(\xi) x_2}.
\end{align}
For simplicity, let the first line in \eqref{eq:hu1} be denoted by $\hu_{11}(x_2;\xi)$ and the second line by $\hu_{12}(x_2;\xi)$. 
As PML is imposed in the region $x_2 > 1$, the closed form of $\hu_1$ in
$x_2\in(0,1)$ is not needed here. In the following, we shall prove that the
truncation error of PML terminating $u_1$ could converge only algebraically.

\subsection{PML-truncated problem}
Now, we introduce a PML in $x_2\in (H,H+L)$ by complexifying $x_2$ via \eqref{eq:sigma}. Recall $\tP = \int_H^{H+L}\alpha(x_2) dx_2$.
The PML-truncated problem is characterized by
\begin{align*}
\partial^2_{x_1} v + \frac{\partial_{x_2}}{1+\bi \sigma(x_2)} \left[\frac{\partial_{x_2}v}{1+\bi\sigma(x_2)}\right]  + (k^2+\epsilon \sin x_1 h(x_2))  v =& -\delta(x - x^*),\quad {\rm in}\quad\mathbb{R}_+^2,\\
v =& 0,\quad{\rm on}\quad \Gamma\cup\Gamma_{H+L}.
\end{align*}
The error function $w(x)=v(x)-u(x_1,\tx_2)$ is governed by
\begin{align*}
\partial^2_{x_1} w + \frac{\partial_{x_2}}{1+\bi \sigma(x_2)} \left[\frac{\partial_{x_2}w}{1+\bi\sigma(x_2)}\right]  + (k^2+\epsilon \sin x_1 h(x_2))  w =& 0,\quad {\rm in}\quad\mathbb{R}_+^2,\\
w(x_1,0) =& 0,\\
w(x_1,H+L) =& -u(x_1,\tP).
\end{align*}
Similar as before, we assume $w = \sum_{j=0}^{\infty} w_j\epsilon^j$ as the problem is well-posed and its solution is analytic at $\epsilon=0$.
It can be seen that the leading term $w_0$ solves
\begin{align}
 \partial^2_{x_1} w_0 + \frac{\partial_{x_2}}{1+\bi \sigma(x_2)} \left[\frac{\partial_{x_2}w_0}{1+\bi\sigma(x_2)}\right]  + k^2w_0 =& 0, x_2>0\\
w_0(x_1,0) =& 0,\\  
w_0(x_1,H+L) =& -u_0(x_1,\tP).   
\end{align}
Using the method of Fourier transform, it can be seen that since $x_2>x_2^*$,
\[
\hw_0(x_2;\xi) =  \bi e^{\bi \mu(\xi) \tP} \frac{\sin(\mu(\xi)x_2^*)}{\mu(\xi)}\frac{\sin(\mu(\xi) \tx_2)}{\sin(\mu(\xi) \tP)}.
\]
Next, $w_1$ is governed by 
\begin{align*}
 \partial^2_{x_1} w_1 + \frac{\partial_{x_2}}{1+\bi \sigma(x_2)} \left[\frac{\partial_{x_2}w_1}{1+\bi\sigma(x_2)}\right]   + k^2w_1 =& -\sin x_1 h(x_2) w_0, x_2>0\\
w_1(x_1,0) =& 0,\\  
w_1(x_1,H+L) =& -u_1(x_1,\tP).   
\end{align*}
Again, we take the Fourier transform of the above equations,
\begin{align}
    \frac{1}{1+\bi \sigma(x_2)}\left[\frac{\hw_1'}{1+\bi\sigma(x_2)}\right]'   + \mu^2(\xi)\hw_1 
    =& \frac{\bi }{2}h(x_2) \left[\hw_0(x_2;\xi+1) - \hw_0(x_2;\xi-1)\right],\\
    \hw_1(0;\xi) =& 0,\\
    \hw_1(H+L;\xi) =& -\hu_1(\tP;\xi).
\end{align}
For $x_2>1$, we assume
\[
\hw_1(x_2;\xi) = A_1(\xi) e^{\bi \mu(\xi) \tx_2} + B_1(\xi) e^{-\bi \mu(\xi) \tx_2},
\]
with the two unknowns $A_1$ and $B_1$.
For $x_2\in(0,1)$, we find
\begin{align*}
\hw_1^-(x_2;\xi) =& -\frac{1}{2}\Bigg[
e^{\bi \mu(\xi+1) \tP} \frac{\sin(\mu(\xi+1)x_2^*)}{\mu(\xi+1)}\frac{\sin(\mu(\xi+1) x_2)}{\sin(\mu(\xi+1) \tP)(2\xi+1)}\\
&+e^{\bi \mu(\xi-1) \tP} \frac{\sin(\mu(\xi-1)x_2^*)}{\mu(\xi-1)}\frac{\sin(\mu(\xi-1) x_2)}{\sin(\mu(\xi-1) \tP)(2\xi-1)} \Bigg]
\end{align*}
is a special solution such that
\[
\hw_1(x_2;\xi) = \hw_1^-(x_2;\xi) + C_1(\xi) \sin(\mu(\xi) x_2).
\]
with the unknown $C_1(\xi)$.
The continuity conditions on $x_2=1$ and the boundary condition on $x_2=L+D$ imply
\begin{align*}
    A_1(\xi) e^{\bi \mu(\xi) \tP} + B_1(\xi) e^{-\bi \mu(\xi)\tP } =& -\hu_1(\tP;\xi),\\
    A_1(\xi) e^{\bi \mu(\xi) } + B_1(\xi) e^{-\bi \mu(\xi) } =& \hw_1^-(1;\xi) + C_1(\xi) \sin(\mu(\xi)),\\
    A_1(\xi) e^{\bi \mu(\xi) } - B_1(\xi) e^{-\bi \mu(\xi) } =& \frac{[\hw_1^-]'(1;\xi)}{\bi\mu(\xi)} - \bi C_1(\xi) \cos(\mu(\xi)).
\end{align*}
On solving the above equation, 
\begin{align*}
C_1(\xi) =& -\frac{\hu_1(\tP;\xi) }{\sin(\mu(\xi)\tP)}\\
&+\frac{\cos(\mu(\xi)[\tP-1])}{2\sin(\mu(\xi)\tP)}\Bigg[
e^{\bi \mu(\xi+1) \tP} \frac{\sin(\mu(\xi+1)x_2^*)}{\mu(\xi+1)}\frac{\sin(\mu(\xi+1) )}{\sin(\mu(\xi+1) \tP)(2\xi+1)}\\
&+e^{\bi \mu(\xi-1) \tP} \frac{\sin(\mu(\xi-1)x_2^*)}{\mu(\xi-1)}\frac{\sin(\mu(\xi-1) )}{\sin(\mu(\xi-1) \tP)(2\xi-1)} \Bigg]\\
&+\frac{\sin(\mu(\xi)[\tP-1])}{2\sin(\mu(\xi)\tP)}\Bigg[
e^{\bi \mu(\xi+1) \tP} \frac{\sin(\mu(\xi+1)x_2^*)}{\mu(\xi)}\frac{\cos(\mu(\xi+1) )}{\sin(\mu(\xi+1) \tP)(2\xi+1)}\\
&+e^{\bi \mu(\xi-1) \tP} \frac{\sin(\mu(\xi-1)x_2^*)}{\mu(\xi)}\frac{\cos(\mu(\xi-1) )}{\sin(\mu(\xi-1) \tP)(2\xi-1)} \Bigg].
\end{align*}
It is left to justify  for any $x_2\in(0,1)$ and a finite fixed $x_1\in\mathbb{R}$, 
\[
w_1(x) = \frac{1}{2\pi}\int_{-\infty}^{+\infty}\hw_1(x_2;\xi) e^{-\bi \xi x_1} d\xi  
\]
decays algebraically when $k$ is a half-integer and exponentially otherwise as $|\tP|\to\infty$.
\subsection{Algebraic convergence}
Taking the inverse Fourier transform of $\hw_1$, we obtain  
\[
w_1(x) = \frac{1}{2\pi}\int_{-\infty}^{+\infty}\hw_1(x_2;\xi)e^{-\bi \xi x_1}d\xi = w_1^A(x) + w_1^E(x),
\]
where we have defined
\begin{align*}
w_1^{A}(x) =& \frac{1}{2\pi}\int_{-\infty}^{+\infty}\frac{-\sin(\mu(\xi)x_2) \hu_1(\tP;\xi) }{\sin(\mu(\xi)\tP)}e^{-\bi \xi x_1}d\xi\\
w_1^{E}(x) =& \frac{1}{2\pi}\int_{-\infty}^{+\infty}\left[\hw_1 + \frac{\sin(\mu(\xi)x_2) \hu_1(\tP;\xi) }{\sin(\mu(\xi)\tP)}\right]e^{-\bi \xi x_1}d\xi.
\end{align*}
We now show that the error function $w_1$ decays algebraically for $k=\frac{1}{2}$ as the branch points of $\mu(\xi+1)$ and $\mu(\xi)$ coincide when $|\xi+\frac{1}{2}|\ll 1$ and the branch points of $\mu(\xi-1)$ and $\mu(\xi)$ coincide when $|\xi-\frac{1}{2}|\ll 1$.
\begin{mytheorem}
    Let $k=\frac{1}{2}$ and $\tP=P(\alpha+\bi\beta)$ for some constants
    $\alpha,\beta>0$, and let $x_1$ be such that $(2\pi)^{-1}x_1\notin \mathbb{Z}$. Then,
    \begin{align}
     w_1(x) = \tP^{-4}\frac{-\pi^3x_2x_2^*}{45}\sin\frac{x_1}{2} + {\cal O}(P^{-5}),\quad P\to+\infty.   
    \end{align}
\begin{proof}
  We justify that $w_1^A$ decays algebraically with $P$. Using the same idea, it
  is straightforward to derive that $w_1^E$ decays exponentially. 
 By \eqref{eq:hu1} and by $\hu_{11}(x_2;\xi)=-\hu_{12}(x_2;-\xi)$,
\begin{align}
    w_1^{A}&=-\frac{1}{2\pi}\int_{-\infty}^{\infty}\frac{\sin(\mu(\xi)x_2)\hu_{11}(\tP;\xi)}{\sin(\mu(\xi)\tP)}e^{-\bi\xi x_1}d\xi{-\frac{1}{2\pi}\int_{-\infty}^{\infty}\frac{\sin(\mu(\xi)x_2)\hu_{12}(\tP;\xi)}{\sin(\mu(\xi)\tP)}e^{-\bi\xi x_1}d\xi}\nonumber\\
    &=-\frac{1}{2\pi}\int_{-\infty}^{\infty}\frac{\sin(\mu(\xi)x_2)\hu_{11}(\tP;\xi)}{\sin(\mu(\xi)\tP)}e^{-\bi\xi x_1}d\xi{+\frac{1}{2\pi}\int_{-\infty}^{\infty}\frac{\sin(\mu(\xi)x_2)\hu_{11}(\tP;\xi)}{\sin(\mu(\xi)\tP)}e^{\bi\xi x_1}d\xi}
\nonumber\\
&=\int_{-\infty}^{+\infty}\frac{\bi}{2\pi} 
f(\xi)e^{\bi\mu(\xi+1)x_2^*+\bi \mu(\xi) \tP}\frac{\sin(\mu(\xi)x_2)\sin(\xi x_1)}{ \sin(\mu(\xi)\tP)} d\xi\nonumber\\
&=\left(\int_{|\xi + 1/2|<\epsilon_0} + \int_{|\xi - 1/2|<\epsilon_0} + \int_{\min|\xi\pm 1/2|>\epsilon_0} \right)\cdots d\xi\nonumber\\
&=:I_1(x;x_2^*) + I_2(x;x_2^*) + I_3(x;x_2^*),
\end{align}
where $\epsilon_0>0$ is a sufficiently small constant such that $\sin(x_1 \xi) >0 $ for any $|\xi+1/2|<\epsilon_0$, and we have defined
\[
f(\xi) = \frac{1}{(2\xi+1)}\left[
\frac{\sin(\mu(\xi+1))}{\mu(\xi+1)}\cos(\mu(\xi)) -\frac{\sin(\mu(\xi))}{\mu(\xi)}\cos(\mu(\xi+1))\right].
\]
It can be easily verified that $f(\xi)$ is an analytic function on $\mathbb{R}$
with $f(-1/2)=-1/3$. It is straightforward to verify that for $|\xi\pm
1/2|>\epsilon$, there exists a positive constant $\delta_0$, that depends only on $\epsilon_0$, such that 
\[
|I_3(x;x_2^*)|\leq C e^{-\delta_0 P}.
\]
As for $I_2(x;x_2^*)$, we deform the path $1/2-\epsilon_0 \to 1/2+\epsilon_0$
to the lower half circle $\{z: |z-1/2|=\epsilon_0, {\rm Im}(z)\leq 0\}$ in
$\mathbb{C}^{+-}$. It is then straightforward to deduce that $I_2$ decays
exponentially with $P$. 

It is left to prove that $I_1(x;x_2^*)$ decays algebraically as $P\to+\infty$.
Note that the idea of path deformation breaks down here as the whole line
segment $(-1/2-\epsilon_0,-1/2+\epsilon_0)$ always overlaps with one of branch cuts of $\mu(\xi+1)$ and
$\mu(\xi)$ \cite{zha22}. Now, let
\begin{align*}
 I_{11}(x;x_2^*) =& \int_{-1/2}^{-1/2+\epsilon_0} \frac{\bi}{2\pi} 
f(\xi)e^{\bi\mu(\xi+1)x_2^*+\bi \mu(\xi) \tP}\frac{\sin(\mu(\xi)x_2)\sin(\xi x_1)}{ \sin(\mu(\xi)\tP)} d\xi,\\   
     I_{12}(x;x_2^*) =& \int_{-1/2-\epsilon_0}^{-1/2} \frac{\bi}{2\pi} 
f(\xi)e^{\bi\mu(\xi+1)x_2^*+\bi \mu(\xi) \tP}\frac{\sin(\mu(\xi)x_2)\sin(\xi x_1)}{ \sin(\mu(\xi)\tP)} d\xi. 
\end{align*}
For $I_{11}$, let $\xi = -1/2+t^2$ such that $t\in(0,\sqrt{\epsilon_0})$. Then,
\begin{align*}
    I_{11}(x;x_2^*) =& \int_{0}^{\sqrt{\epsilon_0}} F_{11}(t) \frac{e^{\bi t\sqrt{1-t^2} \tP}t^2 }{ \sin(t\sqrt{1-t^2}\tP)} dt,
\end{align*}
where we have defined
\[
F_{11}(t) = \frac{\bi}{\pi} 
f(-1/2+t^2)e^{-t\sqrt{1+t^2}x_2^*}\frac{\sin(t\sqrt{1-t^2}x_2)\sin((-1/2+t^2) x_1)}{t},
\]
analytic at $t=0$. We further introduce a new variable $s = t\sqrt{1-t^2}$ with
$s\in (0,\sqrt{\epsilon_0-\epsilon_0^2})$ to transform
\[
I_{11}(x;x_2^*)  = \int_{0}^{\sqrt{\epsilon_0}} G_{11}(s)\frac{s^2 }{ 1-e^{-2\bi s\tP}} ds,
\]
where
\[
G_{11}(s)=\frac{2\bi F_{11}(t(s))t'(s)}{1-t^2(s)},
\]
analytic at $s=0$. We first claim
\[
\left|\int_{0}^{\sqrt{\epsilon_0}} \frac{[G_{11}(s) - G_{11}(0) - G_{11}'(0) s]s^2 }{ 1-e^{-2\bi s\tP}} ds\right|= {\cal O}(P^{-5}),\quad P\to\infty.
\]
This can be seen by 
\begin{align*}
    {\rm L.H.S.} \leq& C\int_{0}^{\sqrt{\epsilon_0}} \frac{s^4 }{ |1-e^{-2\bi s\tP}|} ds
    \leq C\int_{0}^{\sqrt{\epsilon_0}} \frac{s^4 }{ e^{2s P\beta}-1} ds\\
    \leq& CP^{-5}\int_{0}^{+\infty}\frac{s^4}{e^{2s \beta}- 1} ds,
\end{align*}
for some generic constant $C$.
Next, apply the Cauchy integral formula, we get
\begin{align*}
 &\int_{0}^{\sqrt{\epsilon_0}} [G_{11}(0) + G_{11}'(0) s]\frac{s^2 }{ 1-e^{-2\bi s\tP}} ds\\
 =&\tP^{-3}\int_{0}^{\sqrt{\epsilon_0} \tP} G_{11}(0)\frac{s^2 }{ 1-e^{-2\bi s}} ds + \tP^{-4}\int_{0}^{\sqrt{\epsilon_0} \tP} G_{11}'(0)\frac{s^3 }{ 1-e^{-2\bi s}} ds\\
 \sim &\tP^{-3}G_{11}(0)\int_{0}^{(\alpha+\bi\beta)\infty} \frac{s^2 }{ 1-e^{-2\bi s}} ds + \tP^{-4}G_{11}(0)\int_{0}^{(\alpha+\bi\beta)\infty} \frac{s^3 }{ 1-e^{-2\bi s}} ds\\
 =&\tP^{-3}G_{11}(0)\int_{0}^{+\infty} \frac{-\bi s^2 }{ 1-e^{2 s}} ds + \tP^{-4}G_{11}'(0)\int_{0}^{+\infty} \frac{s^3 }{ 1-e^{2 s}} ds,\quad P\to\infty.
\end{align*}
Thus,
\begin{align*}
    I_{11}(x;x_2^*) = \tP^{-3}G_{11}(0) \int_{0}^{+\infty} \frac{-\bi s^2 }{ 1-e^{2 s}} ds + \tP^{-4}G_{11}'(0) \int_{0}^{+\infty} \frac{s^3 }{ 1-e^{2 s}} ds + {\cal O}(P^{-5}),\quad P\to\infty.
\end{align*}
Following the same approach,  we get 
\begin{align*}
    I_{12}(x;x_2^*) = \tP^{-3}G_{12}(0)\int_{0}^{+\infty} \frac{s^2 }{ 1-e^{2 s}} ds + \tP^{-4}G_{12}'(0)\int_{0}^{+\infty} \frac{s^3 }{ 1-e^{2 s}} ds + {\cal O}(P^{-5}),\quad P\to\infty.
\end{align*}
By directly computing the involved coefficients and the equation
\[
\int_{0}^{+\infty} \frac{s^3 }{ 1-e^{2 s}} ds = -\frac{\pi^4}{240},
\]
we get the desired result.
\end{proof}
\end{mytheorem}

As $w_1$ decays algebraically in the physical region $\{x:x_2\in(0,1)\}$ as $\tP\to\infty$, one can make $\epsilon$ sufficiently small to ensure that $\epsilon w_1$ becomes the dominant term of $w$ considering $w_0$ decays exponentially. Consequently, there exists a compact region $D$ such that  
\[
\max_{x\in D} |w(x)| \sim C P^{-4}, \quad{\rm for}\quad P\gg 1,
\]
for some constant $C>0$ depending on $\epsilon$. 
\begin{myremark}
    {For $2k\in\mathbb{Z}^+\backslash\{1\}$, it can be shown in a similar fashion that $w_1$ decays exponentially. Nevertheless, one shall see that there exists some integer $j\geq 2$, such that the high-order error term $w_j$ decays algebraically. In other words, the error function $w$ shall always decay algebraically as long as $k$ is a half-integer.}
\end{myremark}
\section{Numerical examples}
In this section, we carry out several experiments to validate the previously established theory. In all examples, the period of the scattering surface is set to be $2\pi$ the same as before. Thus, the half-integers become the exceptional case where the convergence rate of the PML is downgraded to be algebraic. To observe such a phenomenon, we require the PML truncation error dominates the numerical error so that the accuracy of the numerical solution becomes essential. Thus, the recently developed high-accuracy PML-BIE method \cite{yuhulurat22} becomes a suitable solver to check the phenomenon.  Basically, the PML-BIE method separates $\Omega_H$ into unit cells, establish 
 BIEs on the boundary of the unit cell containing the perturbed part $\Gamma\backslash\Gamma_0$, and then evaluates $\tu$ elsewhere via Green's representation theorem; in Appendix A, we present the basic idea of this numerical solver.

As depicted in Figure~\ref{fig:surfaces},
\begin{figure}[htbp]
\centering
\includegraphics[width=0.4\columnwidth]{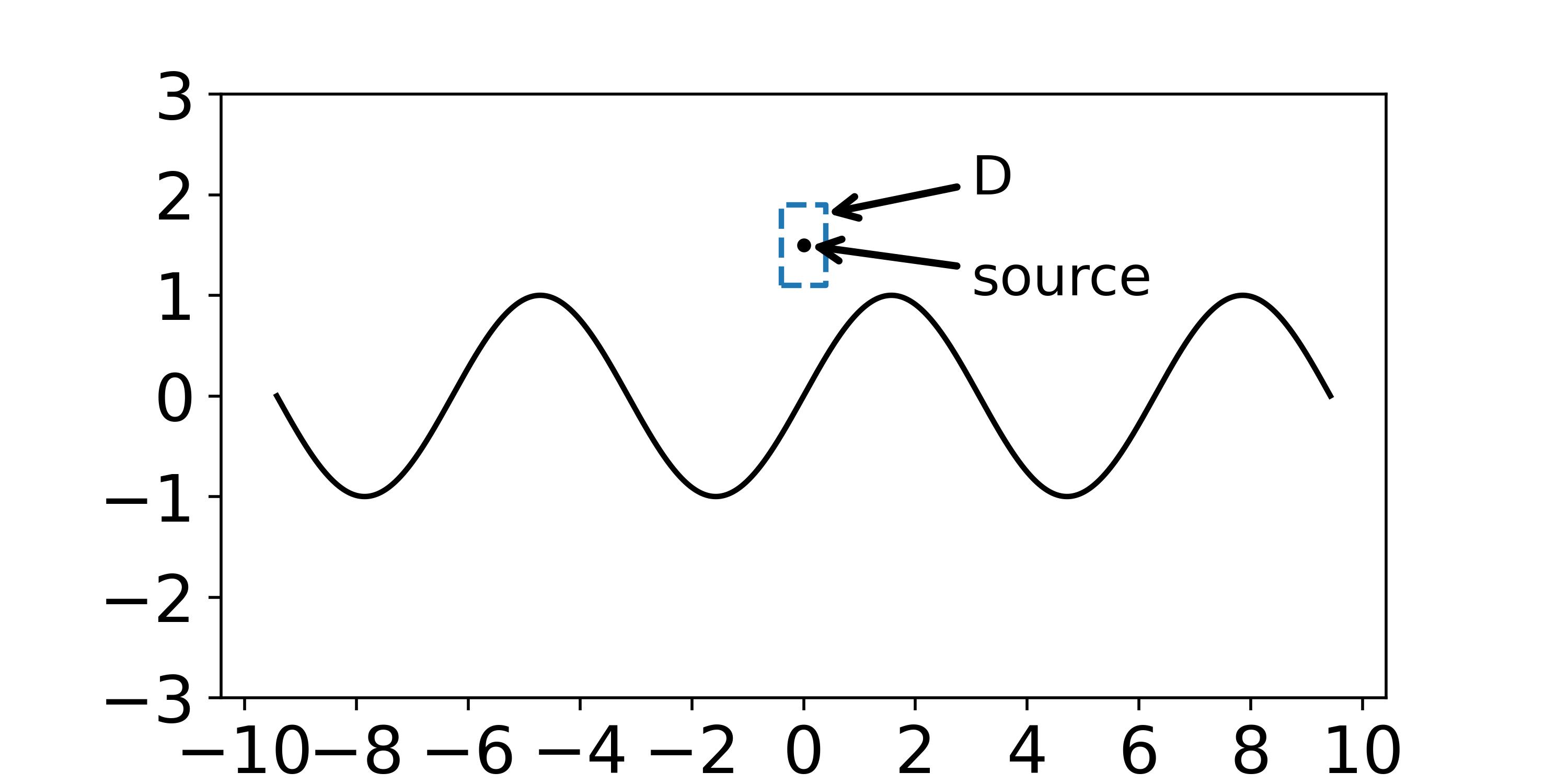}
\includegraphics[width=0.4\columnwidth]{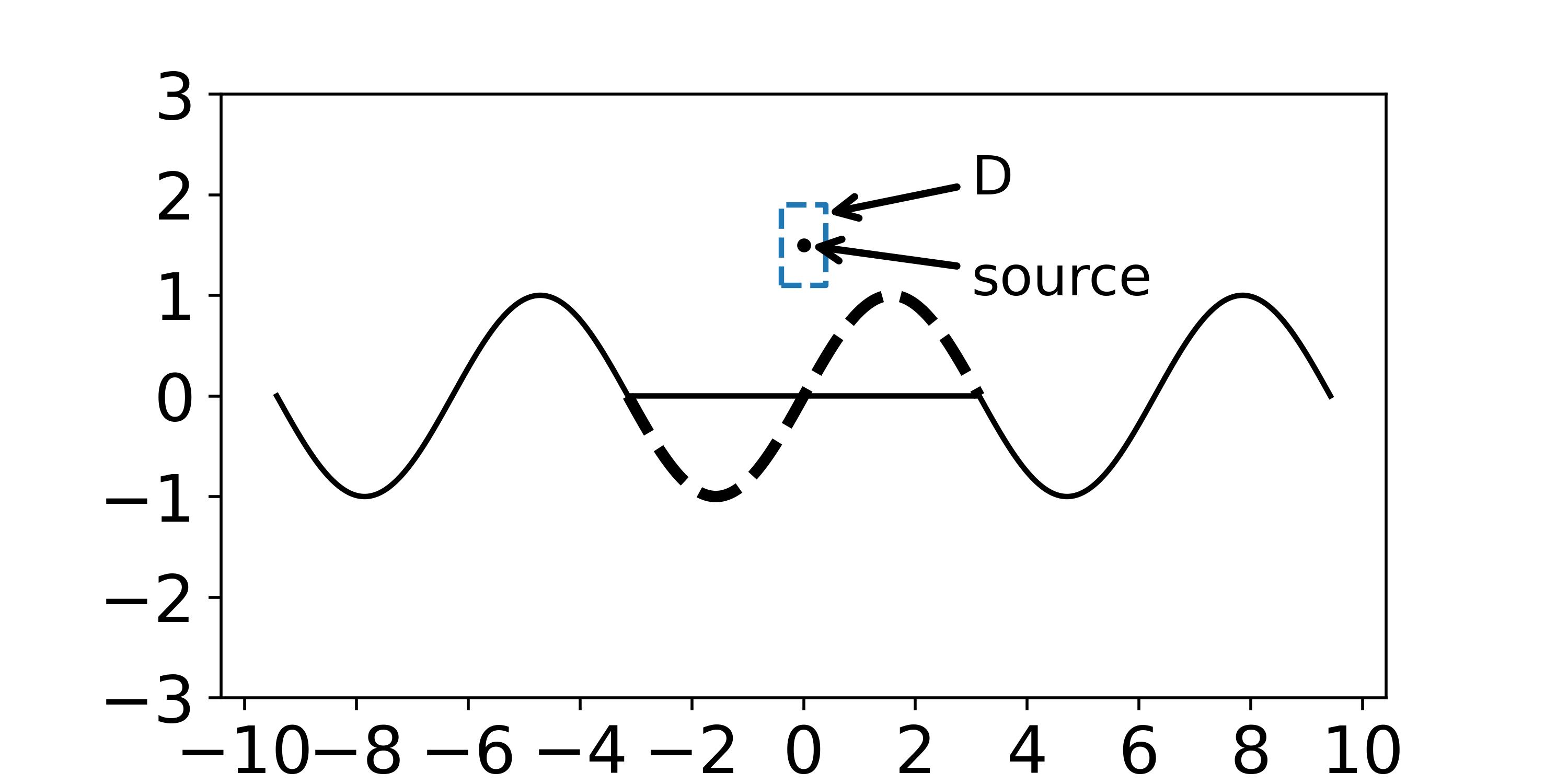}\\
\includegraphics[width=0.4\columnwidth]{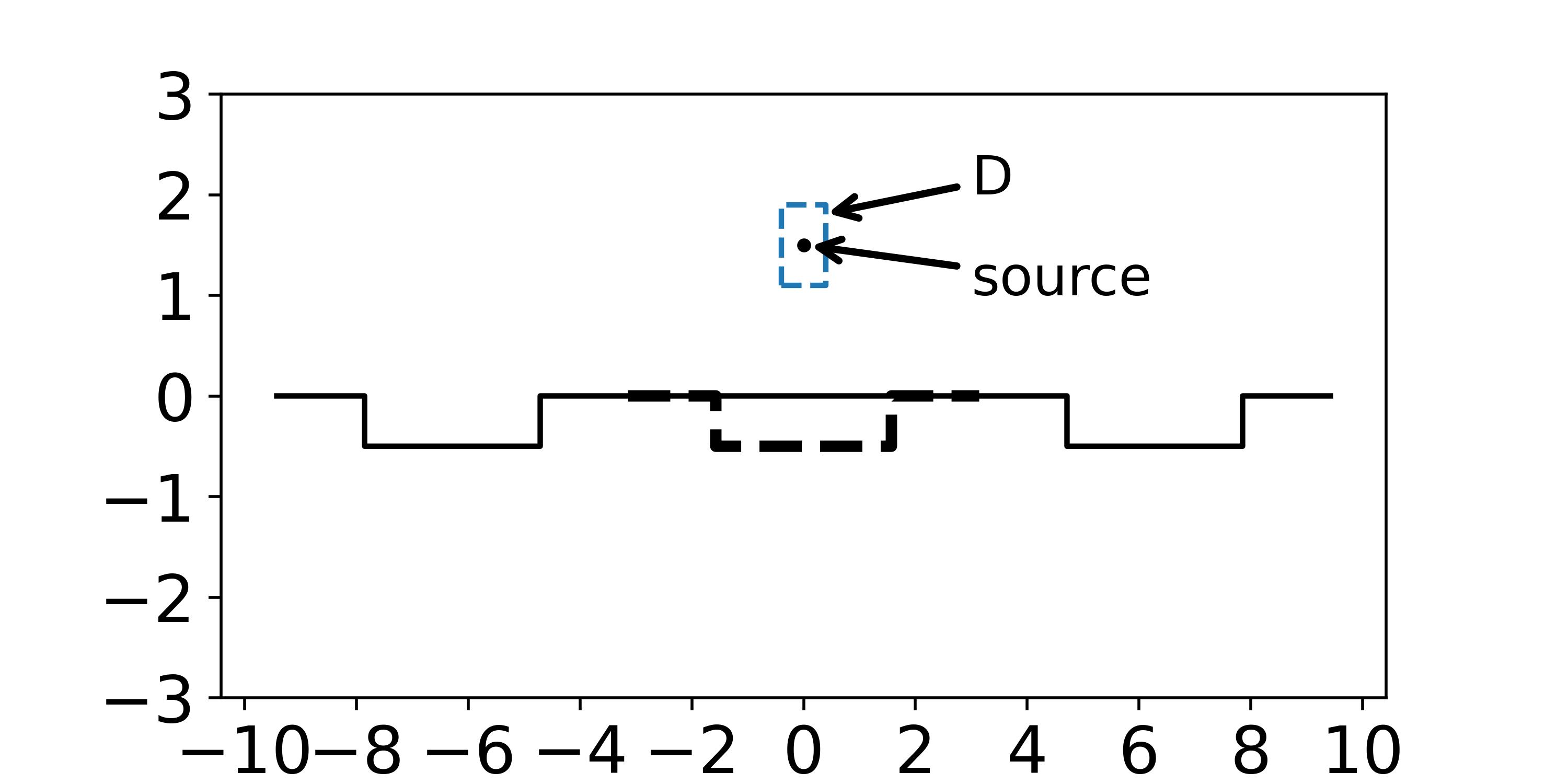}
\includegraphics[width=0.4\columnwidth]{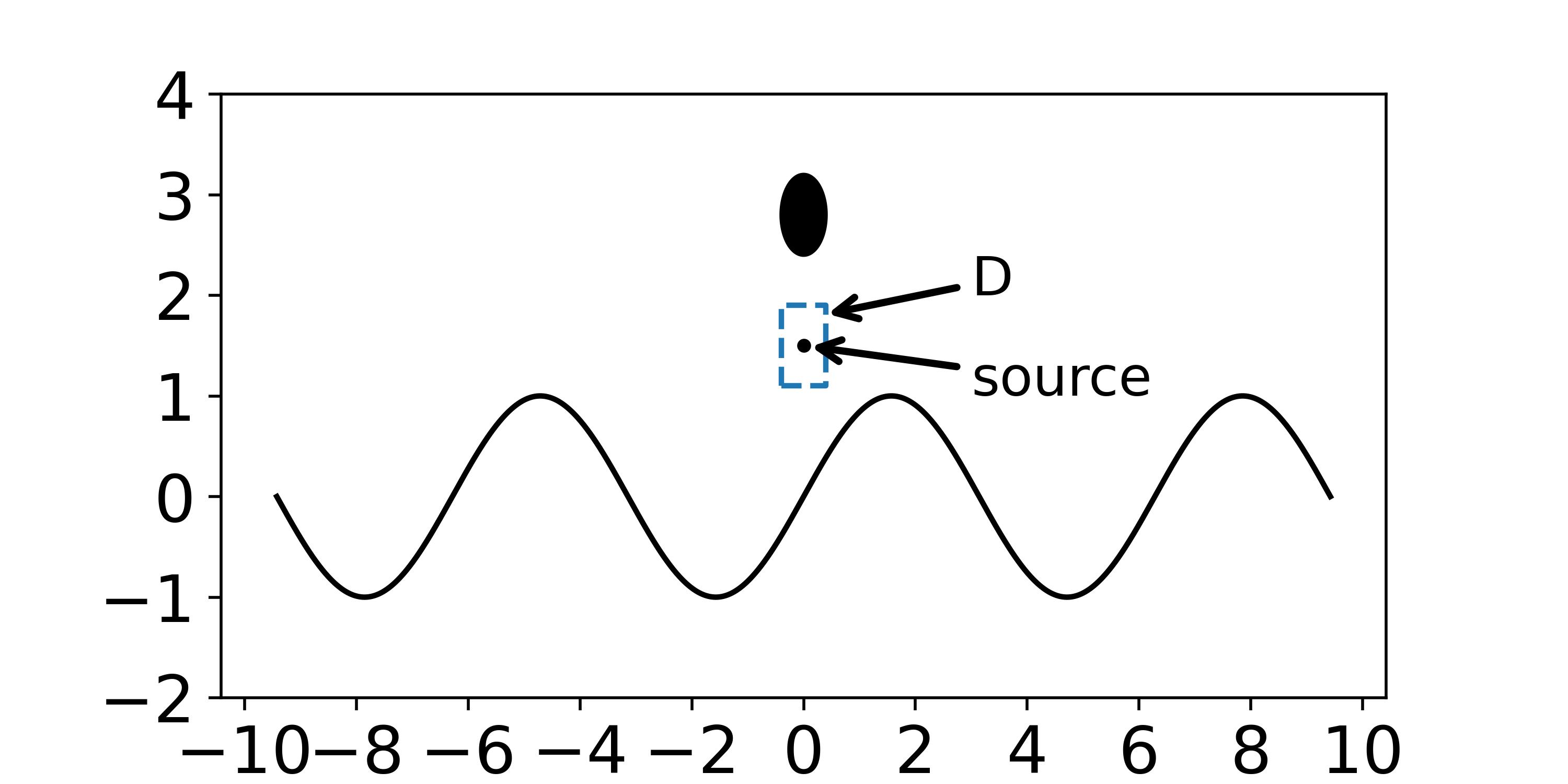}
\caption{Four locally perturbed surfaces: $\Gamma_1$ (top left), $\Gamma_2$ (top right), $\Gamma_3$ (bottom left), and $\Gamma_4$ (bottom right). Dashed lines indicate the perturbed parts; dashed rectangles indicate the regions $D$ where numerical solutions are computed; the dots indicate the locations of the source points; the dark region indicates a  non-penetrable obstacle with zero Dirichlet condition on its boundary.}
\label{fig:surfaces}
\end {figure}
we consider four different surfaces in the following:\\
{
\begin{itemize}
    \item[$\Gamma_1$:]  A sine curve $x_2=\sin x_1$;
    \item[$\Gamma_2$:] A sine curve $x_2=\sin x_1$ locally perturbed by the straight line $\{x:x_2=0,x_1\in(-\pi,\pi)\}$;
    \item[$\Gamma_3$:] A locally perturbed binary grating.
    \item[$\Gamma_4$:] The union of a sine curve $x_2=\sin x_1$ and the boundary of a non-penetrable obstacle occupying the region $x_1^2 + (x_2-2.8)^2\leq 0.4^2$.
\end{itemize}
}
To setup the PML, we take
\[
\sigma(x_2)=\left\{\begin{array}{lc}
\frac{2f^6_2}{f^6_1+f^6_2},&x_2\in[H,H+L/2],\\
2,&x_2\geq H+L/2,
\end{array}\right.
\]
in (\ref{eq:sigma}) to ensure that $\tilde{u}$ is sufficiently smooth across $x_2=H$, where 
\[f_1=(\frac{1}{2}-\frac{1}{6})\xi^3+\frac{\xi}{6}+\frac{1}{2},\ f_2=1-f_1,\ \xi=\frac{2x_2-(2H+L/2)}{L/2}.\]
In all examples, we consider only point-source incidences at the same source point $(0,1.5)$, and compute numerical solutions in $D=[-0.3,0.3]\times[1.2,1.8]$, sufficiently away from the three aforementioned scattering surfaces $\{\Gamma_j\}_{j=1}^3$, to ease the PML-BIE method for accurately computing $u$ in $D$. 
We take $H=3$ for the first three surfaces and $H=4$ for the last surface $\Gamma_4$, and $S=2.8$ in all examples, use a sufficiently refined mesh in the PML-BIE solver, and let $L$ vary to check the accuracy of the PML. 
A reference solution $\tu^{\rm ref}$ is defined as the numerical solution for a sufficiently large $L$, and the $H^1$-error is then defined by 
\[
E_{\rm rel} = \frac{||\tu - \tu^{\rm ref}||_{H^1(D)}}{||\tu^{\rm ref}||_{H^1(D)}}.
\]
Certainly, a discrete $H^1$-norm is used to approximate the continuous norm $||\tu^{\rm ref}||_{H^1(D)}$ as $\tu^{\rm ref}$ is available on a grid in $D$ in the numerical solver.

To illustrate the affection of the wavenumber $k$ on the decaying rate of PML, we consider two groups of values of $k$: 
\begin{itemize}
    \item[(i)] $k\in\{3.1,3.01,3.001,3\}$.
    \item[(ii)] $k\in\{1.5,3,6\}$;
\end{itemize}
Numerical results are shown in Figures~\ref{fig:g1}, \ref{fig:g2}, \ref{fig:g3}, and \ref{fig:g4}.
\begin{figure}[htbp]
\centering
(a)\includegraphics[width=0.43\columnwidth]{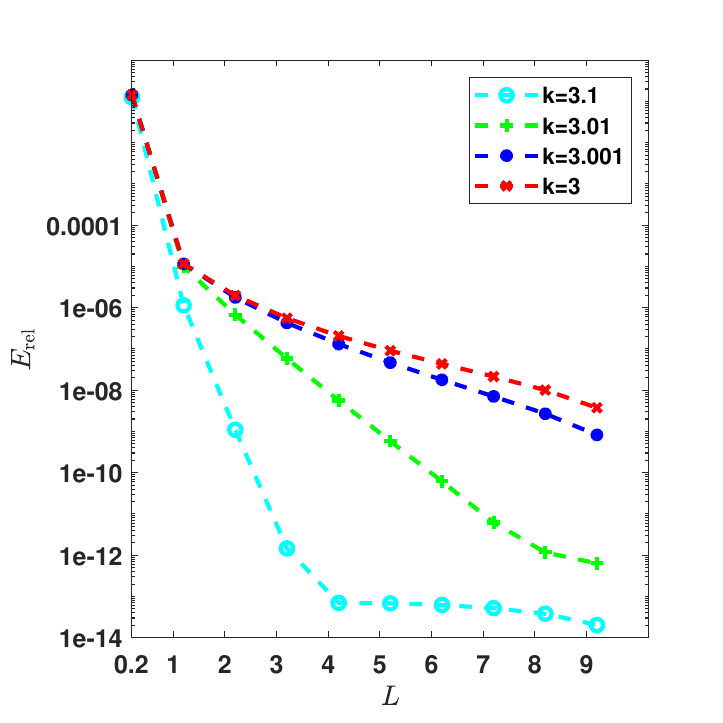}
(b)\includegraphics[width=0.43\columnwidth]{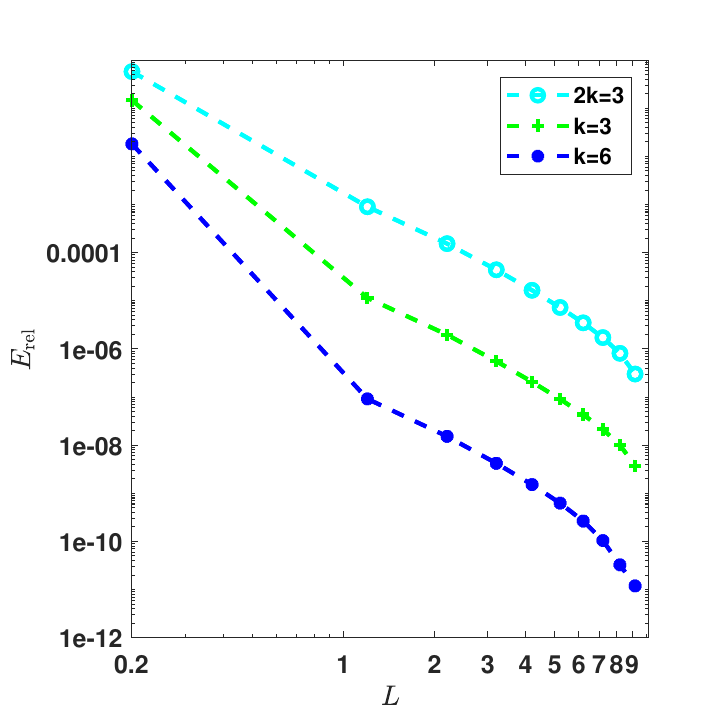}
\caption{Numerical error $E_{\rm rel}$ against PML thickness $L$ for the scattering surface $\Gamma_1$ for different values of $k$. }
\label{fig:g1}
\end {figure}
\begin{figure}[htbp]
\centering
(a)\includegraphics[width=0.43\columnwidth]{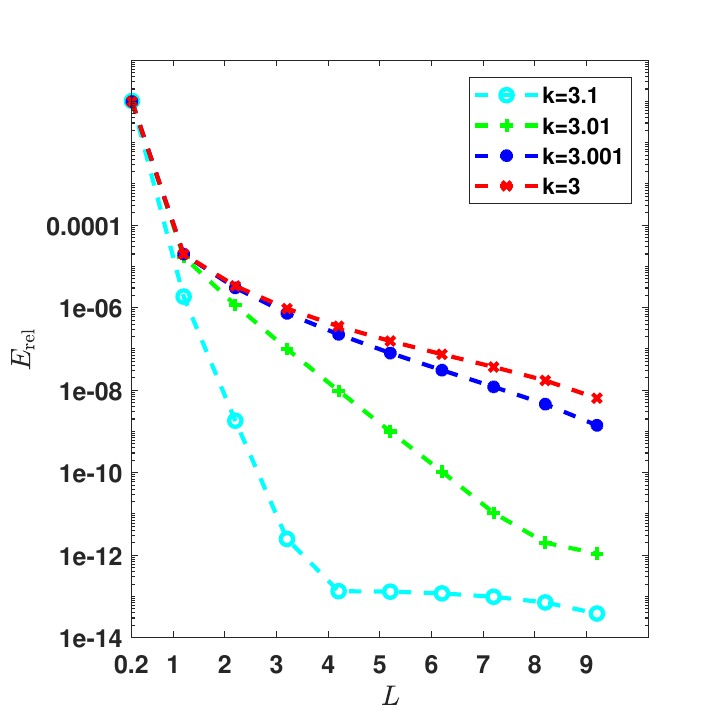}
(b)\includegraphics[width=0.43\columnwidth]{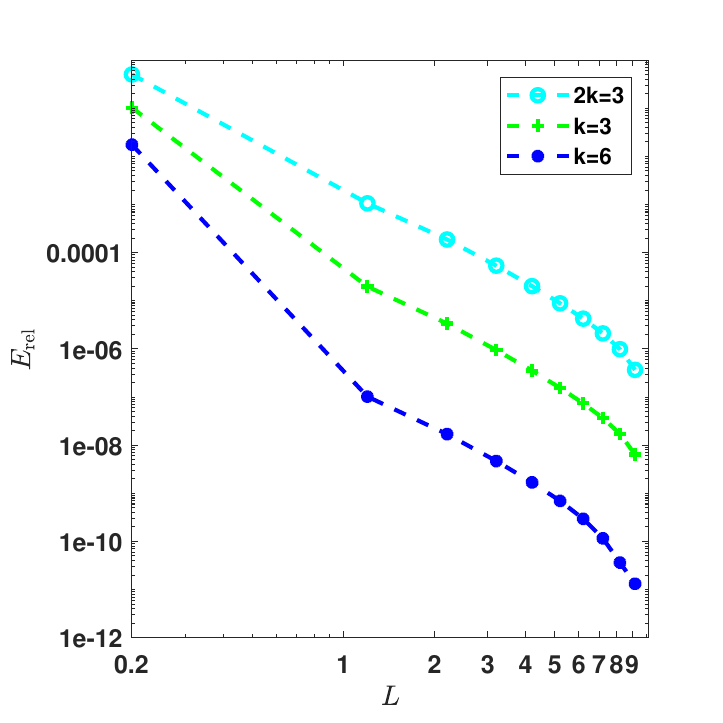}
\caption{Numerical error $E_{\rm rel}$ against PML thickness $L$ for the scattering surface $\Gamma_2$ for different values of $k$. }
\label{fig:g2}
\end {figure}
\begin{figure}[htbp]
\centering
(a)\includegraphics[width=0.43\columnwidth]{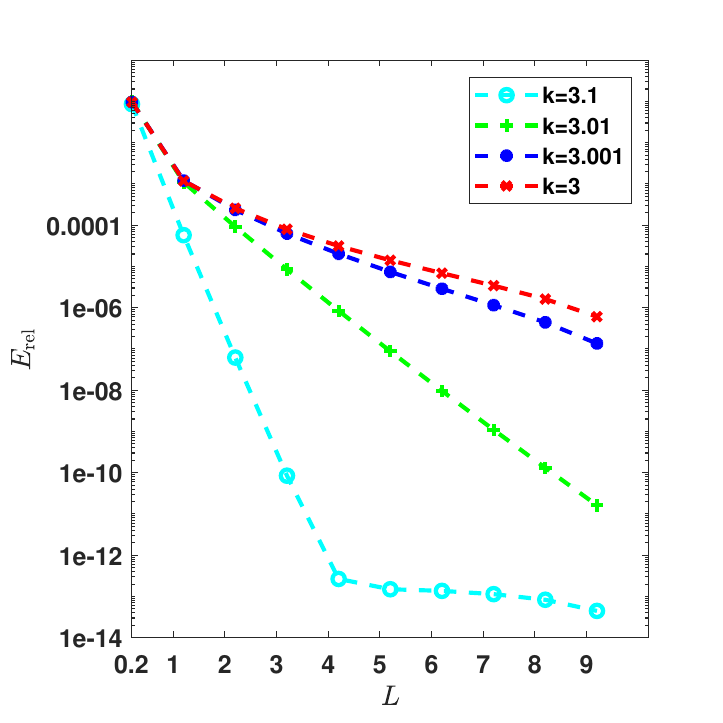}
(b)\includegraphics[width=0.43\columnwidth]{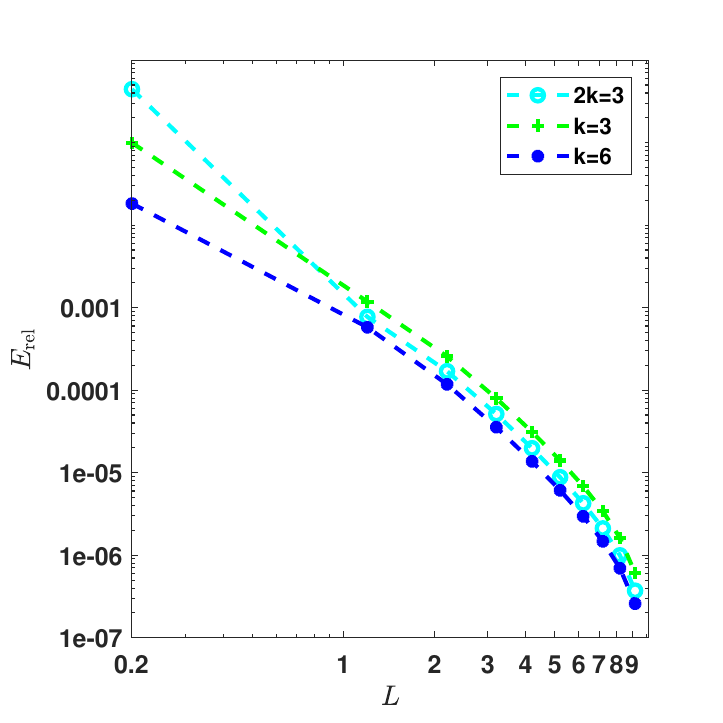}
\caption{Numerical error $E_{\rm rel}$ against PML thickness $L$ for the scattering surface $\Gamma_3$ for different values of $k$. }
\label{fig:g3}
\end {figure}
\begin{figure}[htbp]
\centering
(a)\includegraphics[width=0.43\columnwidth]{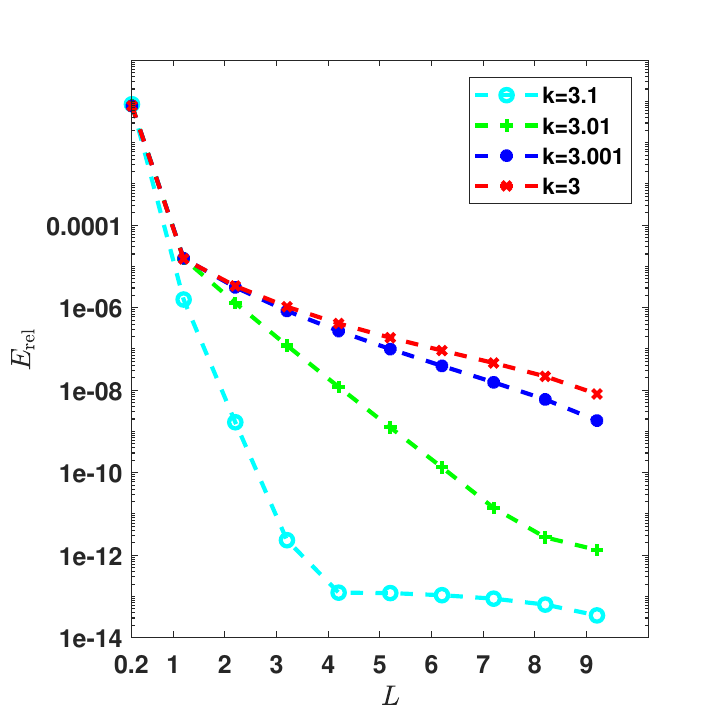}
(b)\includegraphics[width=0.43\columnwidth]{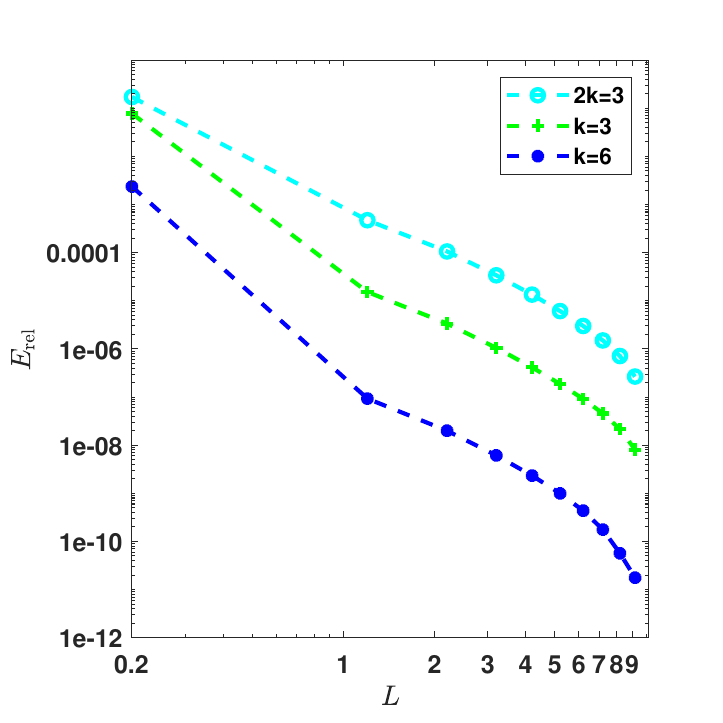}
\caption{Numerical error $E_{\rm rel}$ against PML thickness $L$ for the scattering surface $\Gamma_4$ for different values of $k$. }
\label{fig:g4}
\end {figure}
We make several observations below. Firstly, when $k$ is sufficiently away from half-integers, $E_{\rm rel}$ decays exponentially as the PML thickness $L$ increases. Secondly,  as $k$ approaches $3$, the decaying rate decreases dramatically, making the accuracy goes down from $14$ digits to merely $6$ digits (c.f. Figure~\ref{fig:g3}); this was not observed in the numerical results of \cite{zha22} due to the limited accuracy of the FEM solver. Lastly, the convergence rate seems to be independent of $k$ as $k$ varies in half-integers; heuristically, the decaying exponent is usually proportional to $k$ for PML that is capable of exponentially absorbing outgoing waves. This is a crucial evidence for the algebraically decaying rate for PML truncation errors at half-integer wavenumbers.

\section{Conclusions and Discussions}
This paper established the PML convergence theory for the problem of wave scattering by a locally perturbed periodic surface. For either the original or the PML-truncating problems, we solved an associated scattering problem with the unperturbed periodic surface to construct the Dirichlet-to-Neumann map on a bounded surface that bounds the whole perturbed region. Using the previous PML convergence theory for unperturbed periodic surfaces \cite{zha22}, we justified that the difference between the two DtN maps on the same bounded surface is exponentially small, or algebraically small for half-integer wavenumbers, with the PML parameters. Consequently, the convergence of the PML solution to the true solution in any compact region was established.

We have found the deteriorate of PML in periodic structures as the wavenumber $k$ approaches any half-integer. Our theory indicates that the PML truncation error can at most achieve a fourth-order convergence rate. To ensure the accuracy, PML must be made as thick as possible (c.f. a nine-wavelength thick PML in Figure~\ref{fig:g3} that retrieves $6$ digits only), making itself lose attraction. Thus, a truncation technique that is uniformly accurate for all wavenumbers is desired in practice. We shall investigate this issue in a future work.

\appendix
\section{The PML-BIE method}
In the appendix, we briefly introduce the high-accuracy PML-BIE method developed in \cite{yuhulurat22}. For simplicity, we consider the unperturbed case $\Gamma=\Gamma_0$. For the PML problem (\ref{eq:pmlgov:1}) and (\ref{eq:pmlgov:2}), we assume $f=-\delta(x-x^*)$ for $x^*=(x_1^*,x_2^*)\in\Omega_H$, such that the BIE method is sufficient to get $\tilde{u}$ in $\Omega_H$. As shown in Fig.~\ref{fig:append},
\begin{figure}[htbp]
\centering
(a)\includegraphics[width=0.45\columnwidth]{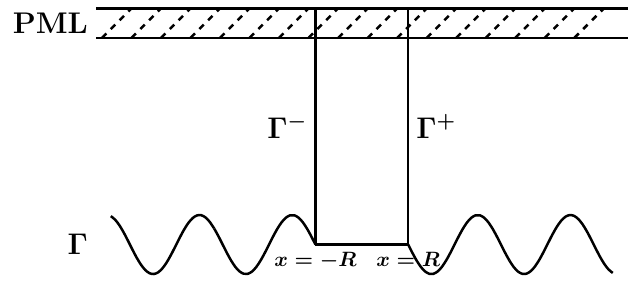}
(b)\includegraphics[width=0.405\columnwidth]{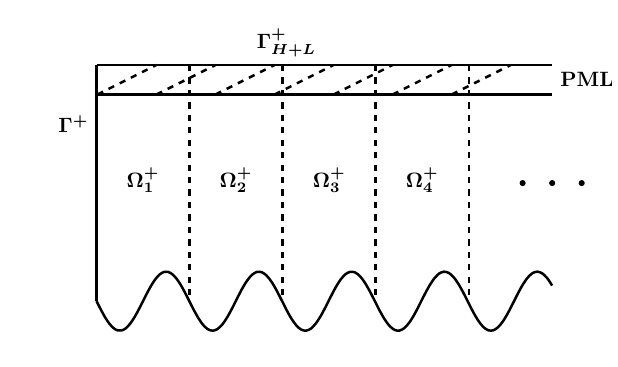}
\caption{Schematic of the periodic structure after the PML truncation.}
\label{fig:append}
\end {figure}
the method basically consists of three steps:
\begin{itemize}
    \item[I.] Divide the domain $\Omega_{H+L}$ into three regions by two vertical lines $x_2=\pm R$ for some sufficiently large $R>0$ with $|x_2^*|<R$; 
    \item[II.] Compute Neumann-to-Dirichlet (NtD) operators ${\cal N}^\pm$ that map $\partial_{\nu} u$ to $\tilde{u}$ on the two boundaries $x_2=\pm R$, where $\nu$ denotes the unit outer normal; 
    \item[III.] Solve the resulting boundary value problem in $\Omega_{H+L}\cap\{x:|x_2|<R\}$.
\end{itemize}
Step II is essential as it truncates the unbounded domain $\Omega_{H+L}$. Without loss of generality, we compute the NtD operator ${\cal N}^+$ on $\Gamma^+$. In doing so, we split the periodic domain $\Omega_{H+L}\cap\{x:x_2>R\}$ into identical unit cells $\Omega_n^+:=\Omega_{H+L}\cap\{x: R+2\pi n<x_2<R+2\pi(n+1)\}, n=0,1,\cdots$. Let $\Gamma_n^+=\Omega_{H+L}\cap \{x:x_2=R+2\pi n\}$.
We define the marching operator ${\cal R}^+: H^{-1/2}(\Gamma^+)\to H^{-1/2}(\Gamma^+)$ (Here, we suppress the subscript of $\Gamma_n^+$ as the related Sobolev spaces are independent of $n$) that maps $\partial_{x_1} u$ on $\Gamma_n^+$ to itself on $\Gamma_{n+1}^+$; it is proved that ${\cal R}^+$ does not depend on $n$ and that $||{\cal R}^+||< 1$. On the other hand, in each unit cell, due to the Dirichlet boundary conditions on $\Gamma_{H+L}$ and $\Gamma$, we find the NtD operators ${\cal N}^{(0)}_{ij}: H^{-1/2}(\Gamma^+)\to \widetilde{H^{1/2}}(\Gamma^+), i,j=1,2$ for one unit cell that satisfy
\begin{equation}
   \label{eq:appd:1} 
   \left[
\begin{matrix}
{\cal N}^{(0)}_{11} & {\cal N}^{(0)}_{12} \\
{\cal N}^{(0)}_{21} & {\cal N}^{(0)}_{22} \\
\end{matrix}
\right]
\left[
\begin{matrix}
\partial_{x_1}\tu^+_n \\
\partial_{x_1}\tu^+_{n+1} \\
\end{matrix}
\right]=\left[
\begin{matrix}
\tu^+_n \\
\tu^+_{n+1} \\
\end{matrix}
\right],
\end{equation}
where $\tu^+_{n}$ is the trace of $\tu$ on $\Gamma_n^+$, and $\partial_{x_1} \tu^+_n$ can be regarded as the normal derivative of $\tu$ on $\Gamma_n^+$. Then, ${\cal R}^+$ is governed by the following Riccati equation
\[
{\cal N}^{(0)}_{21} + {\cal N}^{(0)}_{22}{\cal R}^+ = {\cal N}^{(0)}_{11}{\cal R}^+ + {\cal N}^{(0)}_{12}[{\cal R}^+]^2.
\]
In fact, the above procedure also applies if we use the NtD operator ${\cal N}_{ij}^{(m)}$ for $2^m$ consecutive unit cells; we can iteratively obtain ${\cal N}_{ij}^{(m+1)}$ from ${\cal N}_{ij}^{(m)}$ based on the continuity of $\tu$ and $\partial_{x_1}\tu$ on $\Gamma^+$. Then,  we get
\[
{\cal N}^{(m)}_{21} + {\cal N}^{(m)}_{22}[{\cal R}^+]^{2^m} = {\cal N}^{(m)}_{11}[{\cal R}^+]^{2^m} + {\cal N}^{(m)}_{12}[{\cal R}^+]^{2^{m+1}},
\]
or
\begin{equation}
\label{eq:appd:2}
  [{\cal R}^+]^{2^m} = ({\cal N}^{(m)}_{22}-{\cal N}^{(m)}_{11})^{-1}[{\cal N}^{(m)}_{12}[{\cal R}^+]^{2^{m+1}} - {\cal N}^{(m)}_{21}].  
\end{equation}
Let $M$ be a sufficiently large such that $[{\cal R}^+]^{2^{M+1}}\approx 0$. Eq.~\eqref{eq:appd:2} provides a backward iteration to approximate ${\cal R}^{2^m}, m=M,\cdots, 0$. 
 Then, $\tu_0^+ = [{\cal N}^{(0)}_{11} + {\cal N}^{(0)}_{12} {\cal R}^+]\partial_{x_1}\tu_0^+$ so that 
\[
{\cal N}^+ = {\cal N}^{(0)}_{11} + {\cal N}^{(0)}_{12} {\cal R}^+.
\]
One similarly ${\cal N}^-$ on $x_2=-R$.
Numerically, to approximate ${\cal N}^+$, the most significant step is to approxiamte ${\cal N}_{ij}^{(0)}$ in (\ref{eq:appd:1}). As PML is involved in the unit cell $\Omega_n^+$, the high-accuracy PML-BIE method originated in \cite{luluqia18} is a suitable way to approximate ${\cal N}^{0}$ that maps $\partial_{\nu}\tu$ to $\tu$ on the boundary $\partial \Omega_n^+$. Then, an algebraic manipulation based on the boundary condition on $\Gamma\cup\Gamma_{H+L}$ gets numerical approximations of ${\cal N}_{ij}^{(0)}$ so that ${\cal R}^+$ and ${\cal N}^+$ are approximated. Once we get ${\cal N}^\pm$, the resulting boundary value problem can be solved easily via a standard BIE formulation. Note that $G(x;x^*)=\frac{\bi}{4}H_0^{(1)}(k|x-x^*|)$ must be extracted from the total field $\tu$ first to eliminate the singularity from $f=-\delta(x-x^*)$.
\bibliographystyle{plain}
\bibliography{wt}
\end{document}